\documentclass[11pt,a4paper,leqno,twoside]{amsart}
\usepackage{latexsym,amssymb,amsmath}
\input xy
\xyoption{all}

\def\CC{\mathbb C}

\def\RR{\mathbb R}
\def\HH{\mathbb H}
\def\AA{{\mathbb A}}

\def\OO{\mathbb O}

\def\11{\mathbf 1}
\def\PP{\mathbb P}

\def\e1{\varepsilon_1}
\def\e2{\varepsilon_2}
\def\e3{\varepsilon_3}

\def\P2{{\PP}^2}

\def\00{\underline{0}}
\def\J0{{\cal J}_3(\underline{0})}

\def\PJ0{\PP({\cal J}_3(\underline{0}))}

\def\e{\varepsilon}

\def\AP2{{\AA\PP}^2}
\def\RP2{{\RR\PP}^2}
\def\CP2{{\CC\PP}^2}
\def\HP2{{\HH\PP}^2}
\def\OP2{{\OO\PP}^2}

\newtheorem{theo}{Theorem}[section]
\newtheorem{coro}[theo]{Corollary}
\newtheorem{lemm}[theo]{Lemma}
\newtheorem{prop}[theo]{Proposition}

\theoremstyle{remark}
\newtheorem{rema}[theo]{Remark}

\begin{document}
\title[On K\"{o}the's normality question for central division
algebras]{On K\"{o}the's normality question for locally
finite-dimensional central division algebras}


\author{Ivan D. Chipchakov}
\address{Institute of Mathematics and Informatics\\Bulgarian
Academy of Sciences\\Acad. G. Bonchev Str., Bl. 8\\1113 Sofia,
Bulgaria:\\ E-mail: chipchak@math.bas.bg}

\begin{abstract}
This paper considers K\"{o}the's question of whether every
associative locally finite-dimensional (abbr., LFD) central
division algebra $R$ over a field $K$ is a normally locally finite
(abbr., NLF) algebra over $K$, that is, whether every nonempty
finite subset $Y$ of $R$ is contained in a finite-dimensional
central $K$-subalgebra $\mathcal{R} _{Y}$ of $R$. It shows that
the answer to the posed question is negative if $K$ is a purely
transcendental extension of infinite transcendence degree over an
algebraically closed field $k$. On the other hand, central
division LFD-algebras over $K$ turn out to be NLF in the following
special cases: (i) $K$ is a finitely-generated extension of a
finite or a pseudo-algebraically closed perfect field $K _{0}$;
(ii) $K$ is a higher-dimensional local field with last residue
field equal to $K _{0}$.
\end{abstract}

\maketitle

\par
\medskip
\section{\bf Introduction}

\medskip
Let $F$ be a field, Br$(F)$ the Brauer group of $F$, $s(F)$ the
class of associative finite-dimensional central simple
$F$-algebras, and $d(F)$ the subclass of the division algebras
lying in $s(F)$. It is well-known (cf. \cite[Sects. 12.5 and
14.4]{P}) that Br$(F)$ is an abelian torsion group, so it
decomposes into the direct sum $\oplus _{p \in \mathbb{P}} {\rm
Br}(F) _{p}$ of its $p$-components Br$(F) _{p}$, where
$\mathbb{P}$ is the set of prime numbers. For each $\nabla \in
s(F)$, denote by ind$(\nabla )$ the Schur index of $\nabla $, that
is, the degree deg$(D _{\nabla })$ of the underlying (central)
division $F$-algebra $D _{\nabla }$ of $\nabla $, determined by
Wedderburn's structure theorem (see \cite[Theorem~2.1.6]{He} or
\cite[Sect. 3.5]{P}); also, let exp$(\nabla )$ be the exponent of
$\nabla $, i.e. the order of its (Brauer) equivalence class
$[\nabla ]$, viewed as an element of Br$(F)$. As shown by Brauer
(cf. \cite[Sect. 14.4]{P} or \cite[Sect.~4.5]{GiSz}), exp$(\nabla
)$ divides ind$(\nabla )$ and is divisible by all prime divisors
of ind$(\nabla )$. These relations and the commutativity of
Br$(F)$ are crucial for the proof of Brauer's primary tensor
product decomposition theorem, for every $D \in d(F)$ (cf.
\cite[Proposition~4.5.16]{GiSz} or \cite[Sects. 13.4, 14.4]{P}).
The description of index-exponent pairs over $F$ depends on the
Brauer $p$-dimensions Brd$_{p}(F)$ (in the sense of
\cite[Sect.~4]{ABGV}), $p \in \mathbb{P}$, defined for each $p \in
\mathbb{P}$, as follows: Brd$_{p}(F)$ is finite and equal to
$n(p)$ if $n(p)$ is the least integer $\ge 0$ for which ind$(A
_{p})$ divides exp$(A _{p}) ^{n(p)}$ whenever $A _{p} \in s(F)$
and $[A _{p}] \in {\rm Br}(F) _{p}$; Brd$_{p}(F)$ is infinity (we
write Brd$_{p}(F) = \infty $) if such $n(p)$ does not exist.
Denote by $\Phi _{\rm Br}$ the class of Brauer finite-dimensional
fields (abbr., $\Phi _{\rm Br}$-fields), namely, those fields $E$
for which there are integers $m _{p}(E)\colon p \in \mathbb{P}$,
with $m _{p}(E) \ge {\rm Brd}_{p}(E ^{\prime })$, for every finite
extension $E ^{\prime }/E$ and any $p \in \mathbb{P}$ (compare
with \cite[(2.2) (a)]{Ch6}). When $E \in \Phi _{\rm Br}$ and $p
\neq {\rm char}(E)$, one may take as $m _{p}(E)$ the absolute
Brauer $p$-dimension abrd$_{p}(E)$ of $E$, defined in Section 3
(see Lemma \ref{lemm3.7}). It is known that $\Phi _{\rm Br}$
contains global fields and the field $\mathbb{Q} _{\ell }$ of
$\ell $-adic numbers, for each $\ell \in \mathbb{P}$ (cf.
\cite[(31.4) and (32.19)]{Re}). The problem of finding whether a
field $E _{0}$ belongs to $\Phi _{\rm Br}$ is usually nontrivial
and its consideration lies in the common area of
finite-dimensional division $E _{0}$-algebras and the main topic
of this paper.
\par\smallskip
Let now $A$ be an arbitrary (associative) unital algebra over the
field $F$. We say that $A$ is a locally finite-dimensional (abbr.,
LFD) algebra if its finite subsets generate finite-dimensional
$F$-subalgebras; $A$ is called normally locally finite \cite{Ba}
(abbr., NLF), or totally normal \cite{Sch1}, if every finite
subset $Y \subset A$ is included in a finite-dimensional central
$F$-subalgebra $\Theta _{Y}$ of $A$ that contains the unit of $A$.
The dimension $[A\colon F]$ is said to be countable if it is
countably infinite.
\par\smallskip
Clearly, NLF-algebras over $F$ are LFD and central $F$-algebras.
The early research into infinite-dimensional division NLF-algebras
dates back to K\"{o}the and Schilling (cf. \cite{Ko} and
\cite{Sch1}). K\"{o}the has noted (see \cite[\S 4, page~27]{Ko})
that it is an open question whether central division LFD-algebras
over $F$ are NLF.
\par\smallskip
The present paper gives a negative answer to the stated question
over suitably chosen fields. At the same time, it shows, relying
on the main result of \cite{Ch1}, that central division
LFD-algebras over any Brauer finite-dimensional field are NLF (see
Theorem \ref{theo3.1}). In this respect, our research is motivated
by results proving that a field $F$ is Brauer finite-dimensional
in the following two special cases: (i) $F$ is of type $C _{m}$,
in the sense of Lang \cite{L1}, for some integer $m \ge 0$
(Matzri, see \cite{Mat}); (ii) $F$ is a higher local field, in the
sense of \cite{IHLF}, with a finite last residue field (this
result is mainly due to Khalin, see \cite[Propositions~3.4, 3.5,
and pp.~318-319]{Ch5}, for more details and further references).
In both cases, $F$ is a virtually perfect field, in the sense that
char$(F) = 0$ or char$(F) = q > 0$ and the degree $[F\colon F
^{q}]$ is finite, where $F ^{q} = \{\alpha ^{q}\colon \alpha \in
F\}$.
\par\smallskip
Recall that a field $F$ is said to be of type $C _{m}$ (or a $C
_{m}$-field), for an integer $m \ge 0$, if every $F$-form $f$ (a
homogeneous polynomial $f \neq 0$ with coefficients in $F$) of
degree deg$(f)$ in more than deg$(f) ^{m}$ variables has a
nontrivial zero over $F$. The Diophantine dimension ddim$(F)$ of
$F$ is defined as follows: ddim$(F)$ is the least integer $m \ge
0$, for which $F$ is a $C _{m}$-field; ddim$(F)$ is infinity if
$F$ is not a $C _{m'}$-field, for any $m' \in \mathbb{N}$. For
example, $F$ is a $C _{0}$-field if and only if it is
algebraically closed (cf. \cite[Ch. IX, Theorem~1.4]{L} and
\cite[Theorem~1]{L1}); also, finite fields have type $C _{1}$, by
Chevalley-Warning's theorem (see \cite[Theorem~6.2.6]{GiSz}), and
pseudo algebraically closed (abbr., PAC) perfect fields have type
$C _{2}$ (see \cite{Koll} and \cite[Theorem~21.3.6]{FJ}). It is
known that if $F$ is a $C _{m}$-field and char$(F) = q > 0$, then
$[F\colon F ^{q}] \le q ^{m}$; ddim$(\mathcal{F}) = \infty $ if
$\mathcal{F}$ is a field, char$(\mathcal{F}) = q$ and
$[\mathcal{F}\colon \mathcal{F} ^{q}] = \infty $ (cf.
\cite[Proposition~7.2]{Mat} and \cite[Example~21.2.8]{FJ}). The
class of $C _{m}$-fields is closed under taking algebraic
extensions, for each $m$. When $m \ge 3$ and $F _{0}$ is a perfect
PAC-field, this class contains all extensions of $F _{0}$ of
transcendence degree $m - 2$ (by the Lang-Nagata-Tsen theorem, see
\cite{L1, Na}), and the iterated Laurent (formal power) series
fields $F _{0}((T _{1})) \dots ((T _{m-2}))$ in $m - 2$ variables
over $F _{0}$ (by Greenberg's theorem, see \cite{Gr}).
\par\smallskip
By a $1$-local field, we mean a complete discrete valued field,
and for any integer $m \ge 2$, an $m$-local field with an $m$-th
residue field $K _{0}$ is defined to be a complete field $K _{m}$
with respect to a discrete valuation $w _{0}$, such that the
residue field $\widehat K _{m} := K _{m-1}$ of $(K _{m}, w _{0})$
is an $(m - 1)$-local field with an $(m - 1)$-th residue field $K
_{0}$. If $m \ge 2$ and $v _{m-1}$ is the standard $\mathbb{Z}
^{m-1}$-valued valuation of $K _{m-1}$, then the composite
valuation $v _{m} = v _{m-1} {\ast }w _{0}$ is the standard
$\mathbb{Z} ^{m}$-valued valuation of $K _{m}$. It is known that
$v _{m}$ is Henselian (cf. \cite[Proposition~A.15]{TW}) and $K
_{0}$ equals the residue field of $(K _{m}, v _{m})$. An
informative presentation of the main ideas of the theory of higher
local fields can be found in \cite[Part~I]{IHLF}.
\par
\medskip
\section{\bf The main results}
\par
\medskip
The main purpose of this paper is to prove the following theorem
which gives an answer to K\"{o}the's normality question by showing
that, generally, central division LFD-algebras need not be NLF
(see also Proposition \ref{prop8.2}):
\par
\medskip
\begin{theo}
\label{theo2.1} Let $K _{0}$ be an algebraically closed field, $K
_{\infty }/K _{0}$ a purely transcendental extension of infinite
transcendence degree, $K/K _{\infty }$ a finite field extension, and
$p$ a prime number not equal to {\rm char}$(K _{0})$. Then there
exists a central division {\rm LFD}-algebra $R$ over $K$ with the
following properties:
\par
{\rm (a)} The dimension $[R\colon K]$ is countable, and
finite-dimensional $K$-subalgebras of $R$ are of $p$-power
dimensions;
\par
{\rm (b)} $R$ does not possess finite-dimensional central
$K$-subalgebras different from $K$; in particular, $R$ is not an
{\rm NLF}-algebra over $K$.
\end{theo}
\par
\medskip
\begin{coro}
\label{coro2.2} Let $K _{0}$ be an algebraically closed field, $K
_{0} ^{\prime }/K _{0}$ a purely transcendental extension, and
$K/K _{0} ^{\prime }$ a finite extension. Then the class of
central division {\rm LFD}-algebras over $K$ consists of {\rm
NLF}-algebras if and only if the transcendence degree {\rm
trd}$(K/K _{0})$ is finite.
\end{coro}
\par
\medskip
The left-to-right implication in Corollary \ref{coro2.2} is a
consequence of Theorem \ref{theo2.1}, and the converse implication
is contained in case (a) of the following result:
\par
\medskip
\begin{theo}
\label{theo2.3} Let $K$ be a field and $K _{m}$ an $m$-local field
with a virtually perfect $m$-th residue field $K _{0}$, for some
$m \in \mathbb{N}$. Then central division {\rm LFD}-algebras over
$K$ are {\rm NLF} in the following three cases:
\par
{\rm (a)} {\rm ddim}$(K)$ is finite; in particular, this holds if
$K$ is a finitely-generated extension of a finite or a {\rm PAC}
perfect field.
\par
{\rm (b)} $K/K _{m}$ is an algebraic extension and $K _{0}$ is a
$\Phi _{\rm Br}$-field; this holds, e.g., if {\rm dim}$(K _{0}) <
\infty $ or $K _{0}$ is a finitely-generated extension of a {\rm
PAC}-field.
\par
{\rm (c)} $K/K _{m}$ is a field extension, {\rm trd}$(K/K _{m}) =
1$, and $K _{0}$ is a {\rm PAC} field or an extension of a finite
or algebraically closed field $E$ with {\rm trd}$(K _{0}/E) \le
1$.
\end{theo}
\par
\medskip
The validity of Theorem \ref{theo2.3} (a) in case $K$ is a
finitely-generated extension of a finite field raises interest in
the open question of whether Brauer $p$-dimensions of
finitely-generated extensions of the field $\mathbb{Q}$ of
rational numbers are finite, for every $p \in \mathbb{P}$. In view
of \cite[Theorem~2.1]{Ch4}, a positive answer to this question
will show that, for each pair $d \in \mathbb{N}$, $p \in
\mathbb{P}$, there exists an upper bound $u(p, d) \in \mathbb{N}$
of Brd$_{p}(\Phi _{d})$, for all finitely-generated extensions
$\Phi _{d}$ of $\mathbb{Q}$ with trd$(\Phi _{d}/\mathbb{Q}) = d$.
Therefore, such an answer would prove that the fields $\Phi _{d}$
are contained in $\Phi _{\rm Br}$, which would guarantee the
normality of their central division LFD-algebras.
\par
\medskip
Theorem \ref{theo2.3} is proved in Section 5. Since, by Theorem
\ref{theo3.1} (the main result of Section 3), central division
LFD-algebras over $\Phi _{\rm Br}$-fields are NLF, the idea of our
proof is to show that $K$ is a $\Phi _{\rm Br}$-field in each case
of Theorem \ref{theo2.3}. Section 4 contains preliminaries on
valuation theory, as well as a characterization of the generalized
(Hahn or Mal'cev-Neumann) power series fields lying in $\Phi _{\rm
Br}$, which widens the applicability of Theorems \ref{theo2.3} and
\ref{theo3.1}. Theorem \ref{theo2.3} (b) is included in
Proposition \ref{prop5.1} which describes the intersection of the
class $\Phi _{\rm Br}$ and the one of $m$-local fields. As to
Theorem \ref{theo2.3} (c), it is deduced in Section 5 from the
Harbater-Hartmann-Krashen \cite{HHKr} and Lieblich \cite{Li11}
theorem, an $m$-dimensional generalization of Saltman's theorem
\cite{Salt} on index-exponent relations in central division
algebras over function fields of $\ell $-adic algebraic curves.
Note that the $\Phi _{\rm Br}$-fields singled out by Theorem
\ref{theo2.3} are virtually perfect whereas for each $q \in
\mathbb{P}$, there exists $E _{q} \in \Phi _{\rm Br}$ with char$(E
_{q}) = q$ and $[E _{q}\colon E _{q} ^{q}] = \infty $ (see
\cite[Proposition~2.3 and Remark~8.7]{Ch3}). Therefore, by
\cite[Theorem~2.1~(c)]{Ch4}, Brd$_{q}(E _{q} ^{\prime }) = \infty
$, for every finitely-generated transcendental field extension $E
_{q} ^{\prime }/E _{q}$.
\par
\medskip
To prove Theorem \ref{theo2.1} we first consider the special case
where $K = K _{\infty }$ is a purely transcendental extension of
$K _{0}$ with trd$(K/K _{0})$ countable. This means that $K = K
_{0}(x _{n}, y _{n}\colon n \in \mathbb{N})$, $\{x _{n}, y
_{n}\colon n \in \mathbb{N}\}$ being a set of algebraically
independent elements over $k$. We construct a division $K$-algebra
$R$ as the union $\cup _{n=1} ^{\infty } R _{n}$ of
finite-dimensional division $K _{n}$-algebras $R _{n}$ of
$p$-power dimensions, for every index $n$, where $K _{n} = K
_{0}(x _{j}, y _{j}\colon j = 1, \dots , n)$. The algebras $R
_{n}$, $n \in \mathbb{N}$, are defined inductively so that their
centres $Z _{n}$, $n \in \mathbb{N}$, satisfy the following: (i)
$Z _{n}/K _{0}$ is a purely transcendental extension and trd$(Z
_{n}/K _{0}) = 2n$; (ii) $Z _{n}/K _{n}$ is a Galois extension
with $\mathcal{G}(Z _{n}/K _{n})$ an abelian group of order $p
^{n}$ and period $p$; (iii) $Z _{n} \cap Z _{2n} = K _{n}$ (see
Lemma \ref{lemm8.1}). This implies $R$ is a central division
LFD-algebra over $K$ and Theorem \ref{theo2.1} can be deduced from
our next result.
\par
\medskip
\begin{prop}
\label{prop2.4} Let $K _{0}$ be an algebraically closed field, $p$
a prime number \par\noindent different from {\rm char}$(K _{0})$,
and $\varepsilon _{m}\colon m \in \mathbb{N}$, a sequence of roots
of unity in $K _{0}$, such that $\varepsilon _{1} \neq 1 =
\varepsilon _{1} ^{p}$ and $\varepsilon _{m+1} ^{p} = \varepsilon
_{m}$, for each index $m$. Suppose that $C = K _{0}(x _{i}, y
_{i}\colon i = 1, \dots , n)$ is a purely transcendental extension
of $K _{0}$ with {\rm trd}$(C/K _{0}) = 2n$, for some $n \in
\mathbb{N}$, and $S$ is an algebra over the field \par\noindent $L
= C(\sqrt[p]{y _{i}}\colon i = 1, \dots , n)$, which is isomorphic
to the tensor product $\otimes _{j=1} ^{n} V _{j}$, where $\otimes
= \otimes _{L}$, and for any $j$, $V _{j}$ is the symbol
$L$-algebra $L(X _{j}, \sqrt[p]{Y _{j}}; \varepsilon _{m _{j}})
_{p^{m_{j}}}$ of degree $p ^{m_{j}}$, that is, the $L$-algebra
with generators $\xi _{j}, \eta _{j}$ subject to the relations
$\eta _{j}\xi _{j} = \varepsilon _{m_{j}}\xi _{j}\eta _{j}, \xi
_{j} ^{p^{m_{j}}} = x _{j}, \eta _{j} ^{p^{m_{j}}} = \sqrt[p]{y
_{j}} \in K _{0}$. Then $d(L)$ contains $V _{1}, \dots , V _{n}$
and $S$, $C$ is the only central $C$-subalgebra of $S$, and
\par\vskip0.04truecm\noindent
$[V _{j}\colon L] = p ^{2m_{j}}$, for each $j$.
\end{prop}
\par
\smallskip
Proposition \ref{prop2.4} is proved in Section 7 by methods of
valuation theory. The proof is based on the relations between
valuations and structure of tame division algebras over strictly
Henselian fields (that is, Henselian fields with separably closed
residue fields), considered in Section 6. It is done by showing
(see Proposition \ref{prop7.1}) that $S \otimes _{L} L ^{\prime }
\in d(L ^{\prime })$ and $S \otimes _{L} L ^{\prime }$ does not
possess noncommutative central $C ^{\prime }$-subalgebras,
\par\vskip0.09truecm\noindent
where $C ^{\prime } = K _{0}((X _{1}))((Y _{1})) \dots ((X
_{n}))((Y _{n}))$ and
\par\vskip0.11truecm\noindent
$L ^{\prime } = C ^{\prime }(\sqrt[p]{Y _{i}}\colon i = 1, \dots , n)
= K _{0}((X _{1}))((\sqrt[p]{Y _{1}})) \dots ((X _{n}))((\sqrt[p]{Y
_{n}}))$.
\par\vskip0.11truecm\noindent
The proof of Theorem \ref{theo2.1} in general is presented in
Section 8. It relies on the fact that the above-mentioned
$K$-algebra $R$ can be chosen so that, for any extension $K _{0}
^{\prime }$ of $K _{0}$ linearly disjoint from $K$ over $K _{0}$
(in the sense of \cite[Ch.~VIII, Sect.~3]{L}), $R \otimes _{K} K
^{\prime }$ is a central division {\rm LFD}-algebra over the field
$K ^{\prime } = K _{0} ^{\prime } \otimes _{K _{0}} K$, which does
not possess noncommutative finite-dimensional central $K ^{\prime
}$-subalgebras.
\par
\smallskip The basic notation, terminology and conventions kept
in this paper are standard and essentially the same as in
\cite{TW, L} and \cite{S1}. The notions of an inertial, a nicely
semi-ramified (abbr, NSR), a totally ramified, and a tame
(division) $K$-algebra, where $(K, v)$ is a Henselian field, are
defined in \cite{JW}. Brauer groups and ordered abelian groups are
written additively, Galois groups are viewed as profinite with
respect to the Krull topology, and by a profinite group
homomorphism, we mean a continuous one. Throughout, $\mathbb{Z}$
is the additive group of integers, and $\mathbb{Z} _{\ell }$ is
the additive group of $\ell $-adic integers, where $\ell \in
\mathbb{P}$. For any unital algebra $A$, we consider only
subalgebras of $A$ containing its unit; $Z(A)$, $A ^{\ast }$, and
$A ^{\rm op}$ denote the centre of $A$, the multiplicative group
of $A$, and the algebra opposite to $A$, respectively. When $A \in
s(Z(A))$ and deg$(A) = m$, an extension $Z(A) ^{\prime }$ of the
field $Z(A)$ (cf. \cite[Sect.~12.1]{P}) is called a splitting
field of $A$ if $A \otimes _{Z(A)} Z(A) ^{\prime }$ is isomorphic
as a $Z(A) ^{\prime }$-algebra to the ring $M _{m}(Z(A) ^{\prime
})$ of $m \times m$ matrices with entries in $Z(A) ^{\prime }$.
For any field $E$, $E _{\rm sep}$ is a separable closure of $E$,
$\mathbb{P}_{E} = \{p \in \mathbb{P}\colon p \neq {\rm
char}(E)\}$, and $E ^{\ast n} = \{a ^{n}\colon a \in E ^{\ast
}\}$, for each $n \in \mathbb{N}$. Given a field extension $E
^{\prime }/E$, I$(E ^{\prime }/E)$ stands for the set of
intermediate fields of $E ^{\prime }/E$. In case $E ^{\prime }/E$
is Galois, $\mathcal{G}(E ^{\prime }/E)$ denotes its Galois group;
we say that $E ^{\prime }/E$ is cyclic if $\mathcal{G}(E ^{\prime
}/E)$ is a cyclic group. We write $\mathcal{G}_{E}$ for the
absolute Galois group of $E$, i.e. $\mathcal{G}_{E} =
\mathcal{G}(E _{\rm sep}/E)$, and for any $p \in \mathbb{P}$,
cd$_{p}(\mathcal{G}_{E})$ is the cohomological $p$-dimension of
$\mathcal{G}_{E}$, in the sense of \cite[Ch. I]{S1}, and $E(p)$ is
the maximal $p$-extension of $E$, that is, the compositum of
finite Galois extensions of $E$ in $E _{\rm sep}$ of $p$-power
degrees.
\par
\medskip
\section{\bf Normality of central division LFD-algebras over
$\Phi _{\rm Br}$-fields}
\par
\medskip
The study of central division LFD-algebras over $\Phi _{\rm
Br}$-fields is motivated both by K\"{o}the's normality question
and by the structure theorems for division NLF-algebras over
global fields and $1$-local fields with finite residue fields, due
to Schilling and Barsotti (cf. \cite{Sch1} and \cite{Ba}). These
results, extended in \cite{Ch1}, to all central division
LFD-algebras over such fields, partially generalize Brauer's
(primary tensor product decomposition) theorem. Our main objective
in this section is to demonstrate the applicability of the main
result of \cite{Ch1}, and thereby, to give an affirmative answer
to K\"{o}the's question in the following situation:
\par
\medskip
\begin{theo}
\label{theo3.1} Let $F$ be a $\Phi _{\rm Br}$-field and $R$ a
central division {\rm LFD}-algebra over $F$. Then $R$ is an {\rm
NLF}-algebra and there exist integers $k(p) \ge 0$, $p \in
\mathbb{P}$, and a central $F$-subalgebra $\widetilde R$ of $R$
with the following properties:
\par
{\rm (a)} $\widetilde R$ is $F$-isomorphic to $\otimes _{p \in
\mathbb{P}} R _{p}$, where $\otimes = \otimes _{F}$ and $R _{p}
\in d(F)$ is an $F$-subalgebra of $R$ of degree $p ^{k(p)}$, for
each $p \in \mathbb{P}$;
\par
{\rm (b)} Every finite-dimensional $F$-subalgebra $\mathcal{R}$ of
$R$ is embeddable in $\widetilde R$;
\par
{\rm (c)} For each $p \in \mathbb{P}$, $k(p)$ is the maximal
integer for which there exists $\rho _{p} \in R$ such that $p
^{k(p)} \mid [F(\rho _{p})\colon F]$; also, $\rho _{p}$ can be
chosen to be separable over $F$;
\par
{\rm (d)} $F$ equals the centralizer $C _{R}(\widetilde R) = \{c
\in R\colon c\tilde r = \tilde rc, \tilde r \in \widetilde R\}$.
\par\noindent
Furthermore, if $[R\colon F]$ is countable, then $R$ and
$\widetilde R$ are isomorphic $F$-algebras.
\end{theo}
\par
\medskip
The conclusions of Theorem \ref{theo3.1} are well-known if $R \in
d(F)$, so we tacitly assume throughout its proof that $[R\colon F]
= \infty $. Before proving the theorem, note that, for any
$\lambda \in R$ separable over $F$, there exists a finite subset
$\Sigma _{0}(\lambda )$ of $R$ satisfying the following: $\lambda
\in \Sigma _{0}(\lambda )$ and any finite subset $\Sigma \subset
R$ including $\Sigma _{0}(\lambda )$ generates an $F$-subalgebra
$\Omega _{\Sigma }$ of $R$, such that $[Z(\Omega _{\Sigma
})(\lambda )\colon Z(\Omega _{\Sigma })] = [F(\lambda )\colon F]$.
This statement follows from the finiteness of the set $I(M
_{\lambda }/F)$, where $M _{\lambda }$ is a Galois closure of
$F(\lambda )$ over $F$. When a subset $\Pi \neq \emptyset $ of
$\mathbb{P}$ is finite, and for each $p \in \Pi $, $\rho _{p} \in
R$ is a separable element over $F$ admissible by Theorem
\ref{theo3.1} (c), the statement shows that if $\Sigma _{0}(\rho
_{p}) \subseteq \Sigma $, for all $p \in \Pi $, then $\prod _{p
\in \Pi } p ^{2k(p)}$ divides $[\Omega _{\Sigma }\colon Z(\Omega
_{\Sigma })]$ and $\gcd (\prod _{p \in \Pi } p, [Z(\Omega _{\Sigma
})\colon F]) = 1$. Thus Theorem \ref{theo3.1} (c) makes it
possible to take Step B in the proof of \cite[Proposition~2]{Ch1},
and then to use Brauer's theorem and \cite[Lemma~8.2]{Ch3}, for
proving Theorem \ref{theo3.1} (a).
\par
\medskip
{\it Proof of Theorem 3.1.} We first show that $F$ satisfies the
FC-$p$ condition formulated in \cite{Ch2}, for every $p \in
\mathbb{P}$. Assuming the opposite, one obtains from
\cite[Proposition~2.4]{Ch2}, that there exist an algebraic
extension $F ^{\prime }/F$ and a central division $F ^{\prime
}$-algebra $D ^{\prime } \cong \otimes _{\nu =1} ^{\infty } \Delta
_{\nu }$, which equals the union of $F ^{\prime }$-algebras $D
_{n} = \otimes _{\nu =1} ^{n} \Delta _{\nu }$, $n \in \mathbb{N}$,
where $\otimes = \otimes _{F'}$ and $\Delta _{\nu } \in d(F
^{\prime })$, deg$(\Delta _{\nu }) = p$, for some $p \in
\mathbb{P}$ and all $\nu \in \mathbb{N}$. Hence, by
\cite[(1.3)]{Ch3}, for each $n \in \mathbb{N}$, there is a finite
extension $F _{n}$ of $F$ in $F ^{\prime }$, and an $F
_{n}$-algebra $\Theta _{n} \in d(F _{n})$ with exp$(\Theta _{n}) =
p$, deg$(\Theta _{n}) = p ^{n}$, and $\Theta _{n} \otimes _{F
_{n}} F ^{\prime } \cong D _{n}$ as $F ^{\prime }$-algebras. This
contradicts the assumption that $F \in \Phi _{\rm Br}$ and so
proves that $F$ satisfies conditions FC-$p$, for all $p \in
\mathbb{P}$. Now the existence of integers $k(p) \ge 0$, $p \in
\mathbb{P}$, required by the former part of Theorem \ref{theo3.1}
(c) follows from \cite[Lemma]{Ch1} (or \cite[Lemma~3.9]{Ch2}).
Moreover, for each $p \in \mathbb{P}_{F}$, the latter part of
Theorem \ref{theo3.1} (c) is a consequence of this result and the
well-known fact that $p$ does not divide the degree of any finite
extension of $F$ over its maximal separable subextension (see,
e.g., \cite[Ch.~V, Corollary~6.2]{L}). Suppose now that $p = {\rm
char}(F)$. We show that the latter part of Theorem \ref{theo3.1}
(c) can be deduced by the method of proving the Noether-Jacobson
theorem (cf. \cite[Theorem~3.2.1]{He}); this has been noted
without proof in \cite{Ch1}. Using the method and the double
centralizer theorem (cf. \cite[Theorem~4.3.2]{He}), one obtains
that if $\mathcal{D}$ is a division LFD-algebra with
char$(Z(\mathcal{D})) = p$, and $\lambda _{p} \in (\mathcal{D}
\setminus Z(\mathcal{D}))$ is an inseparable element over
$Z(\mathcal{D})$, then there exists $\xi _{p} \in \mathcal{D}$
separable over $Z(\mathcal{D})$, such that $\xi _{p}\lambda _{p}
^{p} = \lambda _{p} ^{p}\xi _{p}$ and $\lambda _{p}\xi _{p} = (\xi
_{p} + 1)\lambda _{p}$. Therefore, $(\xi _{p} ^{p} - \xi
_{p})\lambda _{p} = \lambda _{p}(\xi _{p} ^{p} - \xi _{p})$, and
$\xi _{p} \notin Z(\mathcal{D})(\xi _{p} ^{p} - \xi _{p})$, so it
follows from the Artin-Schreier theorem (cf. \cite[Sect. VI,
Theorem~6.4]{L}) that $[Z(\mathcal{D})(\xi _{p})\colon
Z(\mathcal{D})(\xi _{p} ^{p} - \xi _{p})] = p$. Thus it turns out
that $p \mid [Z(\mathcal{D})(\xi _{p})\colon Z(\mathcal{D})]$ and
$[Z(\mathcal{D})(\lambda _{p})\colon Z(\mathcal{D})] \mid
[Z(\mathcal{D})(\xi _{p}, \lambda _{p} ^{p})\colon
Z(\mathcal{D})]$. As $\lambda _{p} ^{p^{\mu }}$ is separable over
$Z(\mathcal{D})$, for some $\mu \in \mathbb{N}$ (and $\mu \le
k(p)$), and by \cite[Theorem~4.3.2]{He}, the centralizer $C
_{\mathcal{D}}(\xi _{p})$ lies in the class $d(Z(\mathcal{D})(\xi
_{p}))$, this allows to prove in several similar steps that
$\mathcal{D}$ contains as a $Z(\mathcal{D})$-algebra a finite
separable extension of $Z(\mathcal{D})(\lambda _{p} ^{p^{\mu }})$
of degree divisible by $[Z(\mathcal{D})(\lambda _{p})\colon
Z(\mathcal{D})(\lambda _{p} ^{p ^{\mu }})]$. The obtained result,
applied to the pair $(\mathcal{D}, \lambda _{p}) = (R, \rho
_{p})$, where $\rho _{p} \in R$ and $p ^{k(p)} \mid [F(\rho
_{p})\colon F]$, indicates that $F$ has a finite extension $\Psi
_{p}$ in $F _{\rm sep}$ with $p ^{k(p)} \mid [\Psi _{p}\colon F]$,
which embeds in $R$ as an $F$-subalgebra. This proves Theorem
\ref{theo3.1} (c), since $\Psi _{p}/F$ is simple (cf. \cite[Ch. V,
Theorem~4.6]{L}), and $[F(\tilde \rho _{p})\colon F] = [\Psi
_{p}\colon F]$, for every primitive element $\tilde \rho _{p}$ of
$\Psi _{p}/F$. The rest of our proof relies on the following two
lemmas.
\par
\smallskip
\begin{lemm}
\label{lemm3.2} Let $D$ and $\Theta $ be central division algebras
over a field $E$, such that $D \in d(E)$ and $D \cong \otimes
_{j=1} ^{s} D _{j}$, for some integer $s \ge 2$, where $\otimes =
\otimes _{E}$, and $D _{1}, \dots , D _{s}$ are central
$E$-subalgebras of $D$ of pairwise relatively prime degrees. Then
$D$ embeds in $\Theta $ as an $E$-subalgebra if and only if so do
$D _{1}, \dots , D _{s}$.
\end{lemm}
\par
\smallskip
\begin{proof}
Evidently, the $E$-algebras $D ^{\rm op}$ and $D _{1} ^{\rm op}
\otimes _{E} \dots \otimes _{E} D _{s} ^{\rm op}$ are isomorphic,
so it follows from the Wedderburn-Artin theorem (cf.
\cite[Theorem~2.1.6]{He}), our assumptions, and
\cite[Lemma~3.5~(i), (ii)]{Ch2}, that there exist $m _{0}$ and $m
_{1}, \dots , m _{s} \in \mathbb{N}$ dividing $[D\colon E]$ and
$[D _{1}\colon E], \dots , [D _{s}\colon E]$, respectively, such
that $m _{0} = \prod _{j=1} ^{s} m _{j}$ and $D _{u} ^{\rm op}
\otimes _{E} \Theta \cong M _{m _{u}}(\Theta _{u} ^{\prime })$
over $E$, $u = 0, 1, \dots , s$, for $D _{0} = D$ and some
division $E$-algebras $\Theta _{0} ^{\prime }$, $\Theta _{1}
^{\prime }, \dots , \Theta _{s} ^{\prime }$. Observing also that
$m _{0} = [D\colon E]$ if and only if $m _{j} =[D _{j}\colon E]$,
for $j = 1, \dots , s$, one deduces Lemma \ref{lemm3.2} from
\cite[Lemma~3.5]{Ch2}.
\end{proof}
\par
\medskip
\begin{lemm}
\label{lemm3.3} Let $R$ be a central division {\rm LFD}-algebra
over a field $E$, and let $D \in d(E)$ be a division $E$-algebra
such that $\gcd \{{\rm deg}(D), [E(\alpha )\colon E]\} = 1$, for
each $\alpha \in R$. Then $D \otimes _{K} R$ is a central division
{\rm LFD}-algebra over $E$.
\end{lemm}
\par
\smallskip
\begin{proof}
It follows from our assumptions that $D \otimes _{E} R$ is a
central simple $E$-algebra (cf. \cite[Theorem~4.1.1]{He}), every
finite-dimensional $E$-subalgebra $R _{0}$ of $R$ is a division
algebra, and $\gcd \{[D\colon E], [R _{0}\colon E]\} = 1$ (each
prime divisor $\omega $ of $[R _{0}\colon E]$ divides $[E(r
_{\omega })\colon E]$, for some $r _{\omega } \in R _{0}$, see
\cite[Sect.~13.1]{P} and \cite[Ch.~V, Theorem~4.6 and
Corollary~6.2]{L}, whence $\omega \nmid [D\colon E]$). Since $D
\in d(E)$, this enables one to obtain from \cite[Lemma~3.5
(iii)]{Ch2} (or propositions in \cite[Sects.~13.4 and 14.4]{P})
that $D \otimes _{E} R _{0} \in d(Z(R _{0}))$, which proves Lemma
\ref{lemm3.3}.
\end{proof}
\par
\medskip
We continue with the proof of Theorem \ref{theo3.1}. Theorem
\ref{theo3.1} (c) and \cite[Proposition~2]{Ch1} show that $R$
possesses $F$-subalgebras $R _{p} \in d(F)$, $p \in \mathbb{P}$,
such
\par\noindent
that deg$(R _{p}) = p ^{k(p)}$ and $p \nmid [F(r _{p})\colon F]$,
for any $r _{p} \in C _{R}(R _{p})$ and
\par\noindent
$p \in \mathbb{P}$. Therefore, it can be deduced from Lemma
\ref{lemm3.2} and the Skolem-Noether theorem (cf.
\cite[Theorem~4.3.1]{He}) that $R$ has $F$-subalgebras $T _{n}$,
$n \in \mathbb{N}$, satisfying the following conditions, for each
$n$: $T _{n} \cong \otimes _{j=1} ^{n} R _{p _{j}}$ and $T _{n}
\subseteq T _{n+1}$; $\gcd \{\prod _{\nu =1} ^{n} p _{\nu }, [F(t
_{n}) \colon F]\} = 1$, for every $t _{n} \in C _{R}(T _{n})$.
Here $\otimes = \otimes _{F}$ and $\mathbb{P}$ is presented as the
growing sequence $p _{n}\colon n \in \mathbb{N}$. Clearly, the
union $\widetilde R = \cup _{n=1} ^{\infty } T _{n} := \otimes
_{n=1} ^{\infty } R _{p _{n}}$ has the properties claimed by
Theorem \ref{theo3.1} (a). Further, by the double centralizer
theorem (and \cite[Theorem~4.4.2]{He}), $R = T _{n} \otimes _{F} C
_{R}(T _{n})$ and $Z(C _{R}(T _{n})) = F$, for every $n \in
\mathbb{N}$, which allows to obtain from \cite[Lemma~3.5]{Ch2}
that any $F$-subalgebra $T$ of $R$ with $[T\colon F] < \infty $
embeds in $\widetilde R$ (in fact, $T$ embeds in $T _{n}$ if $p
_{n'} \nmid [T\colon K]$, for any $n' > n$). Thus Theorem
\ref{theo3.1} (b) is proved and the normality of $R$ as an
$F$-algebra reduces to a consequence of Skolem-Noether's theorem;
also, it becomes obvious that $F = \cap _{n=1} ^{\infty } C _{R}(T
_{n}) = C _{R}(\widetilde R)$, as claimed by Theorem
\ref{theo3.1}~(d).
\par\vskip0.04truecm
Suppose finally that $[R\colon F]$ is countable. Since $R$ is NLF
over $F$, it follows from K\"{o}the's theorem (cf. \cite[\S
4]{Ko}) and Theorem \ref{theo3.1} (c) that $R \cong \otimes _{p
\in \mathbb{P}} R _{p} ^{\prime }$, where $\otimes = \otimes
_{F}$, and for each $p \in \mathbb{P}$, $R _{p} ^{\prime }$ is an
$F$-subalgebra of $R$, $R _{p} ^{\prime } \in d(F)$ and deg$(R
_{p} ^{\prime }) \mid p ^{k(p)}$. Denote by $\mathcal{R}_{p}$ the
underlying division $F$-algebra of
\par\vskip0.04truecm\noindent
$R _{p} ^{\prime } \otimes _{F} R _{p} ^{\rm op}$, for each $p$.
Observing that $R = R _{p} \otimes _{F} C _{R}(R _{p}) = R _{p}
^{\prime } \otimes _{F} C _{R}(R _{p} ^{\prime })$
\par\vskip0.04truecm\noindent
and $p \nmid [F(y _{p})\colon F]$, for any $y _{p} \in (C _{R}(R
_{p}) \cup C _{R}(R _{p} ^{\prime }))$, one obtains from the
Wedderburn-Artin theorem, Lemma \ref{lemm3.3} and
\cite[Sect.~12.4, Proposition~b (iv)]{P} that $\mathcal{R}_{p}
\otimes _{F} C _{R}(R _{p} ^{\prime }) \cong C _{R}(R _{p})$ as
$F$-algebras. This yields $\mathcal{R}_{p} = F$ and
\par\vskip0.04truecm\noindent
$R _{p} ^{\prime } \cong R _{p}$ over $F$, for every $p \in
\mathbb{P}$, which proves the concluding assertion of Theorem
\ref{theo3.1}.
\par
\medskip
\begin{rema}
\label{rema3.4} Let $A$ be a simple (unital) algebra over a field
$F$. Then the property of being a simple NLF-algebra (the same as
a locally finite central simple algebra, in the sense of
\cite{BGMV}) is defined by the stronger condition that every
finite subset $Y \subset A$ is included in an $F$-subalgebra
$\Theta _{Y}$ of $A$ with $\Theta _{Y} \in s(F)$. The condition is
satisfied by locally matrix $F$-algebras (the only simple
NLF-algebras if $F$ is algebraically closed). Initiated by
K\"{o}the \cite{Ko}, the research in this area has been continued
by Kurosh, Kurochkin and other authors (see \cite{Ku1}). More
recent results on simple NLF-algebras can be found, for example,
in \cite{BO}, as well as in \cite{BGMV} which gives a
comprehensive presentation of the theory. It is known that central
simple LFD-algebras over an arbitrary field $F$ need not be NLF;
in contrast to results like Corollary \ref{coro2.2}, $F$ always
admits central simple LFD-algebras without noncommutative
$F$-subalgebras lying in $s(F)$ (cf. \cite{Ku}).
\end{rema}
\par
\medskip
Theorem \ref{theo3.1} reduces the research into central division
LFD-algebras over $\Phi _{\rm Br}$-fields to the study of their
finite-dimensional subalgebras. Specifically, the problem of
classifying central division LFD-algebras of countable dimension
over a $\Phi _{\rm Br}$-field $F$, up-to $F$-isomorphisms, is
solved by considering the corresponding problem for algebras $D
_{p} \in d(F)$ of $p$-power degrees, where $p$ runs across
$\mathbb{P}$ (see \cite[Theorem~4.2]{Ch2}, for a complete solution
when $F$ is a global field). The reduction is facilitated by the
known fact (cf. \cite[Sect. 14.4, Proposition~b~(viii)]{P}) that
$d(F)$ contains every tensor product of finitely many algebras
from $d(F)$ of pairwise relatively prime degrees. The fact itself
is a special case of the following result.
\par
\medskip
\begin{prop}
\label{prop3.5} Let $R\{p\}\colon p \in \mathbb{P}$, be central
division {\rm LFD}-algebras over a field $K$, such that $[K(r
_{p})\colon K] = p ^{\nu (r_{p})}$, $\nu (r_{p}) \ge 0$ being an
integer, for any $p \in \mathbb{P}$, $r _{p} \in R\{p\}$. Then the
$K$-algebra $R = \otimes _{p \in \mathbb{P}} R\{p\}$ is central
division and {\rm LFD}.
\end{prop}
\par
\smallskip
\begin{proof}
The assertion holds if and only if, for any finite subset $\Pi =
\{p _{1}, \dots , p _{n}\}$ of $\mathbb{P}$, the $K$-subalgebra
$\Sigma = \otimes _{\nu =1} ^{n} \Sigma _{\nu }$ of $R$ is a
division algebra with
\par\noindent
$Z(\Sigma ) = \otimes _{\nu =1} ^{n} Z(\Sigma _{\nu })$ whenever
$\Sigma _{1}, \dots , \Sigma _{n}$ are finite-dimensional
$K$-subalgebras of $R\{p _{1}\}, \dots , R\{p _{n}\}$,
respectively. Therefore, one proves Proposition \ref{prop3.5},
using the next lemma and arguing by induction on $n$.
\end{proof}
\par
\smallskip
\begin{lemm}
\label{lemm3.6} Let $K$ be a field and $D _{1}$, $D _{2}$ be
finite-dimensional division $K$-algebras, such that $\gcd \{[D
_{1}\colon K], [D _{2}\colon K]\} = 1$. Then $D _{1} \otimes _{K}
D _{2}$ is a division $K$-algebra and $Z(D _{1}) \otimes _{K} Z(D
_{2}) = Z(D _{1} \otimes _{K} D _{2})$.
\end{lemm}
\par
\smallskip
\begin{proof}
Put $Z = Z _{1} \otimes _{K} Z _{2}$, where $Z _{i} = Z(D _{i})$,
$i = 1, 2$. It is known (see \cite[Sect.~12.1, Proposition;
Sect.~15.3, Lemma~a~(i)]{P}) that $Z/K$ is a field extension with
$[Z\colon K] = [Z _{1}\colon K].[Z _{2}\colon K]$. This implies
$[Z\colon Z _{i}] = [Z _{2-i+1}\colon K]$ and $\gcd \{[D
_{i}\colon Z _{i}], [Z\colon Z _{i}]\} = 1$, $i = 1, 2$, so it
follows from \cite[Lemma~3.5]{Ch2} that the $Z$-algebras $D _{1}
\otimes _{Z _{1}} Z := D _{1} ^{\prime }$, $D _{2} \otimes _{Z
_{2}} Z := D _{2} ^{\prime }$ and $D _{1} ^{\prime } \otimes _{Z}
D _{2} ^{\prime }$ lie in $d(Z)$. Observing finally that $D _{1}
^{\prime } \otimes _{Z} D _{2} ^{\prime } \cong D _{1} \otimes
_{K} D _{2}$ as $K$-algebras (cf. \cite[Sect. 9.2, Proposition~
c]{P}), one completes the proof of the lemma.
\end{proof}
\par
\medskip
\label{9899} As noted in Section 2, the results of the next two
sections characterize the $\Phi _{\rm Br}$-fields within the class
of generalized power series fields, and within the class of
$m$-local fields, for some $m \in \mathbb{N}$. Our presentation of
these results relies on the fact that if $E$ is a $\Phi _{\rm
Br}$-field, then its absolute Brauer $p$-dimension
abrd$_{p}(E)$\footnote{The Brauer $p$-dimension, in the sense of
\cite{PS} and \cite{BH}, means the same as the absolute Brauer
$p$-dimension in the present paper.} (defined to be the supremum
of Brd$_{p}(E_{1})$, when $E_{1}$ runs across the set of finite
extensions of $E$ in $E _{\rm sep}$) is finite, for every $p \in
\mathbb{P}$. We refer the reader to \cite{Ch6}, for more results
on the pairs Brd$_{p}(E), {\rm abrd}_{p}(E)$, $p \in \mathbb{P}$,
and to \cite{Ch4}, for their application to the study of
index-exponent relations over finitely-generated transcendental
extensions. The question of whether $E$ is a $\Phi _{\rm
Br}$-field when char$(E) > 0$ and abrd$_{p}(E) < \infty $, for all
$p \in \mathbb{P}$, seems to be open. Our next lemma gives a
positive answer under the condition that $E$ is virtually perfect;
in general, it shows that if $p \in \mathbb{P}_{E}$, then
abrd$_{p}(E)$ equals the supremum of Brd$_{p}(E ^{\prime })$,
where $E ^{\prime }$ runs across the class of finite extensions of
$E$.
\par
\medskip
\begin{lemm}
\label{lemm3.7} Let $E$ be a field of characteristic $q > 0$,
$\overline E$ an algebraic closure of $E _{\rm sep}$, and $E _{\rm
ins}$ the maximal purely inseparable extension of $E$ in
$\overline E$. Then:
\par
{\rm (a)} For each $p \in \mathbb{P}_{E}$, {\rm abrd}$_{p}(E) \ge
{\rm Brd}_{p}(E ^{\prime })$, for every finite extension $E
^{\prime }$ of $E$ in $\overline E$; in addition, {\rm
abrd}$_{p}(E _{\rm ins}) = {\rm abrd}_{p}(E)$;
\par
{\rm (b)} If $E$ is virtually perfect and $[E\colon E ^{q}] = q
^{\delta }$, then {\rm Brd}$_{q}(E ^{\prime }) \le \delta $, for
every finite extension $E ^{\prime }/E$;
\par
{\rm (c)} $E _{\rm ins}$ is a Brauer finite-dimensional field,
provided that so is $E$; the converse holds if $E$ is virtually
perfect.
\end{lemm}
\par
\smallskip
\begin{proof}
Let $E ^{\prime }$ be any finite extension of $E$ in $\overline
E$. Then $E ^{\prime } _{\rm ins}$ equals the compositum $E
^{\prime }.E _{\rm ins}$, and $E ^{\prime }/E _{0} ^{\prime }$ is
a purely inseparable extension, where $E _{0} ^{\prime } = E
^{\prime } \cap E _{\rm sep}$; in particular, $[E ^{\prime }\colon
E _{0} ^{\prime }] = q ^{\nu }$, for some integer $\nu \ge 0$ (cf.
\cite[Ch.~V, Proposition~6.6 and Corollary~6.2]{L}). Therefore, it
follows from Albert-Hochschild's theorem (cf. \cite[Ch. II,
2.2]{S1}) and \cite[Sect.~13.4, Proposition~(vi)]{P} that the
scalar extension map of $s(E _{0} ^{\prime })$ into $s(E ^{\prime
})$ induces an index-preserving group isomorphism Br$(E _{0}
^{\prime }) _{p} \cong {\rm Br}(E ^{\prime }) _{p}$, for each $p
\in \mathbb{P}_{E}$. As finite extensions of $E ^{\prime }$ in $E
^{\prime } _{\rm ins}$ are purely inseparable, their degrees are
$q$-powers, so the scalar extension map $s(E ^{\prime }) \to s(E
^{\prime } _{\rm ins})$ induces an index-preserving isomorphism
Br$(E ^{\prime }) _{p} \cong {\rm Br}(E ^{\prime } _{\rm ins})
_{p}$ as well. These facts, and the one that each finite extension
of $E _{\rm ins}$ in $\overline E$ equals $\widetilde E _{\rm
ins}$, for a suitably chosen finite extension $\widetilde E$ of
$E$ in $E _{\rm sep}$ (cf. \cite[Ch.~V, Theorem~4.6 and Sect.
6]{L}), prove Lemma \ref{lemm3.7} (a). The former part of Lemma
\ref{lemm3.7} (c) is implied by Lemma \ref{lemm3.7} (a), since
every $Y \in I(\overline E/E _{\rm ins})$ is a perfect field,
which ensures that Br$(Y) _{q} = \{0\}$ (see \cite[Ch. VII,
Theorem~22]{A1}); the latter one follows from Lemma \ref{lemm3.7}
(a)-(b). It remains to prove Lemma \ref{lemm3.7} (b), assuming
that $[E\colon E ^{q}] = q ^{\delta }$, for some $\delta \in
\mathbb{N}$. Then $[E ^{\prime }\colon E ^{\prime q}] = q ^{\delta
}$ (see, e.g., \cite[Lemma~2.12]{BH}), and obviously, $E ^{\prime
} _{\rm ins} = \cup _{n=1} ^{\infty } E ^{\prime } _{q^{-n}}$,
where $E ^{\prime } _{q^{-n}} = \{\lambda _{n} \in \overline
E\colon \lambda _{n} ^{q^{n}} \in E ^{\prime }\}$, for each $n$;
also, it is easily obtained that $E ^{\prime } _{q^{-n}}/E
^{\prime }$ is a field extension, $[E ^{\prime }_{q^{-n}}\colon E
^{\prime }] = q ^{\delta n}$, and $E ^{\prime } _{q^{-n}} \subset
E ^{\prime } _{q^{-n-1}}$. At the same time, by \cite[Ch. VII,
Theorem~28]{A1}, $E ^{\prime } _{q^{-n-1}}$ is a splitting field
of any $D _{n} \in d(E ^{\prime } _{q^{-n}})$ with exp$(D _{n}) =
q$. This implies $E ^{\prime } _{q^{-n}}$ is a splitting field of
any $D _{n} ^{\prime } \in d(E ^{\prime })$ with exp$(D _{n}
^{\prime }) = q ^{n}$, proving that deg$(D _{n} ^{\prime }) \mid q
^{\delta n}$ (cf. \cite[Sect.~13.4]{P}) and Brd$_{q}(E ^{\prime })
\le \delta $, as required.
\end{proof}
\par
\medskip
\section{\bf Characterization of the Brauer finite-dimensional fields
which lie in the class of maximally complete equicharacteristic
fields}
\par
\medskip
For any field $K$ with a (nontrivial) Krull valuation $v$, $O
_{v}(K)$ denotes the valuation ring $\{a \in K\colon \ v(a) \ge
0\}$ of $(K, v)$, $M _{v}(K) = \{\mu \in K\colon \ v(\mu ) > 0\}$
the maximal ideal of $O _{v}(K)$, $O _{v}(K) ^{\ast } = \{u \in
K\colon \ v(u) = 0\}$ the multiplicative group of $O _{v}(K)$,
$v(K)$ the value group and $\widehat K = O _{v}(K)/M _{v}(K)$ the
residue field of $(K, v)$, respectively. As usual, $v(K)$ is
assumed to be an ordered abelian group; also, $\overline {v(K)}$
stands for a fixed divisible hull of $v(K)$. The valuation $v$ is
said to be Henselian if it extends uniquely, up-to equivalence, to
a valuation $v _{L}$ on each algebraic extension $L$ of $K$. When
this holds, $(K, v)$ is called a Henselian field. We say that $(K,
v)$ is strictly Henselian if it is Henselian and $\widehat K _{\rm
sep} = \widehat K$. It is well-known that $(K, v)$ is Henselian in
the following two cases: (i) $v(K)$ is embeddable as an ordered
subgroup in the additive group $\mathbb{R}$ of real numbers, and
$K$ is complete with respect to the topology induced by $v$ (cf.
\cite[Ch. XII, Proposition~2.5]{L} or \cite[Theorem~18.3.1]{E3});
(ii) $(K, v)$ is maximally complete, i.e. it does not admit a
valued extension $(K ^{\prime }, v')$ such that $K ^{\prime } \neq
K$, $\widehat K ^{\prime } = \widehat K$ and $v'(K ^{\prime }) =
v(K)$ (cf. \cite[Theorem~15.3.5]{E3}). \label{max_comp} For
example, complete discrete valued fields are maximally complete
(see \cite[Example~3.11]{TW}), and for each $n \in \mathbb{N}$, so
is the $n$-fold iterated Laurent series field $K _{n} = K _{0}((X
_{1})) \dots ((X _{n}))$ over a field $K _{0}$ (that is, $K _{m} =
K _{m-1}((X _{m}))$, for every index $m > 0$) with respect to its
standard valuation $v _{n}$ inducing on $K _{0}$ the trivial
valuation (cf. \cite[Theorem~18.4.1]{E3}). Here $v _{n}(K _{n}) =
\mathbb{Z} ^{n}$, $K _{0}$ is the residue field of $(K _{n}, v
_{n})$, and the abelian group $\mathbb{Z} ^{n}$ is considered with
its inverse-lexicographic ordering (see \cite[Examples~4.2.2 and
9.2.2]{E3}). The condition that $v$ is Henselian has the following
two equivalent forms (cf. \cite[Sect. 18.1]{E3}, or
\cite[Theorem~A.12]{TW}):
\par
\medskip\noindent
(4.1) (a) Given a polynomial $f(X) \in O _{v}(K) [X]$ and an
element $a \in O _{v}(K)$, such that $2v(f ^{\prime }(a)) <
v(f(a))$, where $f ^{\prime }$ is the formal derivative of $f$,
there is a zero $c \in O _{v}(K)$ of $f$ satisfying the equality
$v(c - a) = v(f(a)/f ^{\prime }(a))$;
\par
(b) For each normal extension $\Omega /K$, $v ^{\prime }(\tau (\mu
)) = v ^{\prime }(\mu )$ whenever $\mu \in \Omega $, $v ^{\prime
}$ is a valuation of $\Omega $ extending $v$, and $\tau $ is a
$K$-automorphism of $\Omega $.
\par
\medskip\noindent
Statement (4.1) (a), applied to the polynomials $X ^{\nu } - 1$,
yields the following:
\par
\medskip\noindent
(4.2) If $(K, v)$ is a Henselian field, $\nu \ge 2$ is an integer
and char$(\widehat K) \nmid \nu $, then:
\par
(a) The set $\nabla _{0}(K) = \{\lambda \in K\colon v(\lambda - 1) >
0\}$ is a subgroup of $K ^{\ast \nu }$. An element $a \in O _{v}(K)
^{\ast }$ lies in $K ^{\ast \nu }$ if and only if its residue class
$\hat a$ lies in $\widehat K ^{\nu }$.
\par
(b) If $(K, v)$ is strictly Henselian, then $K$ and $\widehat K$
contain primitive $\nu $-th roots of unity. Moreover, $\widehat K
^{\ast } = \widehat K ^{\ast \nu }$, whence, $K ^{\ast \nu } =
\{\theta _{\nu } \in K ^{\ast }\colon v(\theta _{\nu }) \in \nu
v(K)\}$, which proves the well-known fact that $v$ induces
canonically a group isomorphism $K ^{\ast }/K ^{\ast \nu } \cong
v(K)/\nu v(K)$ (cf. \cite[Lemma~7.78]{TW}, for more details).
\par
\medskip
When $v$ is Henselian, so is $v _{L}$, for any algebraic field
extension $L/K$. In this case, we denote by $\widehat L$ the
residue field of $(L, v _{L})$, and put $v(L) = v _{L}(L)$,
\par\noindent
$O _{v}(L) = O _{v _{L}}(L)$, $M _{v}(L) = M _{v_{L}}(L)$; we
write $v$ instead of $v _{L}$ when there is no danger of
ambiguity. Clearly, $\widehat L/\widehat K$ is an algebraic
extension and $v(K)$ is an ordered subgroup of $v(L)$, such that
$v(L)/v(K)$ is a torsion group; therefore, one may assume without
loss of generality that \label{k99} $v(L)$ is an ordered subgroup
of $\overline {v(K)}$. By Ostrowski's theorem (cf.
\cite[Theorem~17.2.1]{E3}), if $[L\colon K]$ is finite and
$e(L/K)$ is the index of $v(K)$ in $v(L)$, then $[\widehat L\colon
\widehat K]e(L/K)$ divides $[L\colon K]$, and in case $[L\colon K]
\neq [\widehat L\colon \widehat K]e(L/K)$, we have $[L\colon K] =
[\widehat L\colon \widehat K]e(L/K)q ^{d(L/K)}$, for some $d(L/K)
\in \mathbb{N}$, where $q = {\rm char}(\widehat K) > 0$. The
extension $L/K$ is called defectless if $[L\colon K] = [\widehat
L\colon \widehat K]e(L/K)$; it is called totally ramified if
$[L\colon K] = e(L/K)$. We state as a lemma some well-known
criteria for defectlessness of $L/K$:
\par
\medskip
\begin{lemm}
\label{lemm4.1} Let $(K, v)$ be a Henselian field and $L/K$ a
finite extension. Then $[L\colon K] = [\widehat L\colon \widehat
K]e(L/K)$ in the following cases:
\par
{\rm (a)} If char$(\widehat K) \nmid [L\colon K]$ {\rm (}apply
Ostrowski's theorem{\rm )};
\par
{\rm (b)} If $v$ is discrete and $L/K$ is separable {\rm (}see
\cite[Sect. 17.4]{E3}{\rm )};
\par
{\rm (c)} When $(K, v)$ is maximally complete (cf.
\cite[Theorem~31.21]{Wa}).
\par\noindent
In case {\rm (c)}, if {\rm char}$(K) = q > 0$ and $v _{q}$ is the
valuation of $K ^{q}$ induced by $v$, then $(K ^{q}, v _{q})$ is
maximally complete with $v _{q}(K ^{q}) = qv(K)$ and a residue
field $\widehat K ^{q}$. Moreover, $[K\colon K^{q}]$ is finite if
and only if so are $[\widehat K\colon \widehat K ^{q}]$ and $e(K/K
^{q})$.
\end{lemm}
\par
\medskip
The Henselian property of $(K, v)$ ensures (by Schilling's
theorem, see \cite[page~30 and Corollary~1.7]{TW}) that $v$
extends on each $\Delta \in d(K)$ to a unique, up-to equivalence,
valuation $v _{\Delta }$ with the value group $v(\Delta )$. It is
known that $v(\Delta )$ is an ordered abelian group, $v(K)$ is an
ordered subgroup of $v(\Delta )$ of finite index $e(\Delta /K)$,
and the residue division ring $\widehat {\Delta }$ of $(\Delta , v
_{\Delta })$ is a $\widehat K$-algebra. Therefore, one may assume
further, following \cite{TW}, that $v(\Delta )$ is an ordered
subgroup of $\overline {v(K)}$. Note here that, by
Ostrowski-Draxl's theorem \cite{Dr1}, $[\widehat \Delta \colon
\widehat K]e(\Delta /K)$ divides $[\Delta \colon K]$, and when
$[\widehat \Delta \colon \widehat K]e(\Delta /K) \neq [\Delta
\colon K]$, we have char$(\widehat K) = q > 0$ and $[\Delta \colon
K] = q ^{y}.[\widehat \Delta \colon \widehat K]e(\Delta /K)$, for
some $y \in \mathbb{N}$. The $K$-algebra $\Delta $ is called
defectless if $[\Delta \colon K] = [\widehat \Delta \colon
\widehat K]e(\Delta /K)$; we say that it is totally ramified if
$[\Delta \colon K] = e(\Delta /K)$. Lemma \ref{lemm4.1} and the
next lemma give criteria for defectlessness of $\Delta /K$.
\par
\medskip
\begin{lemm}
\label{lemm4.2} Let $(K, v)$ be a Henselian field. Then an algebra
$\Delta \in d(K)$ is defectless over $K$ in the following cases:
\par
{\rm (a)} {\rm char}$(\widehat K)$ $\nmid $ {\rm deg}$(\Delta )$
{\rm (}by the Ostrowski-Draxl theorem{\rm )};
\par
{\rm (b)} $(K, v)$ is maximally complete {\rm (}cf.
\cite[Theorem~3.1]{TY}{\rm )};
\par
{\rm (c)} $v$ is discrete {\rm (}see \cite[Proposition~2.2]{TY} or
\cite[Proposition~4.21 (iii)]{TW}{\rm )}.
\end{lemm}
\par
\medskip
Let $(K, v)$ be a Henselian field and $L/K$ a finite extension. We
say that $L/K$ is tame if it is defectless, char$(\widehat K)
\nmid e(L/K)$, and $\widehat L$ is separable over $\widehat K$;
$L/K$ is said to be inertial if it is tame and $[L\colon K] =
[\widehat L\colon \widehat K]$. Inertial extensions of $K$ are
separable and the following lemma (for its proof, see
\cite[Theorem~A.23 and Corollary~A.25]{TW}) presents some of their
frequently used properties:
\par
\medskip
\begin{lemm}
\label{lemm4.3} Assume that $(K, v)$ is a Henselian field and $K
_{\rm ur}$ is the compositum of inertial extensions of $K$ in $K
_{\rm sep}$. Then:
\par
{\rm (a)} An inertial extension $R ^{\prime }/K$ is Galois if and
only if so is $\widehat R ^{\prime }/\widehat K$. When this holds,
$\mathcal{G}(R ^{\prime }/K)$ and $\mathcal{G}(\widehat R ^{\prime
}/\widehat K)$ are canonically isomorphic.
\par
{\rm (b)} $v(K _{\rm ur}) = v(K)$, $K _{\rm ur}/K$ is a Galois
extension and $\mathcal{G}(K _{\rm ur}/K) \cong
\mathcal{G}_{\widehat K}$.
\par
{\rm (c)} Finite extensions of $K$ in $K _{\rm ur}$ are inertial, and
the natural mapping of $I(K _{\rm ur}/K)$ into $I(\widehat K _{\rm
sep}/\widehat K)$, by the rule $L \to \widehat L$, is bijective.
\par
{\rm (d)} For each finite extension $K _{1}$ of $K$ in $K _{\rm
sep}$, the field $K _{0} = K _{1} \cap K _{\rm ur}$ equals the
maximal inertial extension of $K$ in $K _{1}$; in addition,
$\widehat K _{0} = \widehat K _{1}$.
\end{lemm}
\par
\medskip
Given a field $K _{0}$ and an ordered (nontrivial) abelian group
$\Gamma $, the generalized power series field $K _{0}((\Gamma ))$
is defined below, following Hahn, Mal'cev and Neumann (see
\cite[Sect. 2.8]{E3} or \cite[1.1.4]{TW}): \label{999}
\par
\medskip\noindent
{\bf Definition.} {\it The field $K _{0}((\Gamma ))$ is determined
by $K _{0}$ and $\Gamma $ as the set of those mappings $f\colon K
_{0} \to \Gamma $, whose support supp$(f) := \{\gamma \in \Gamma
\colon f(\gamma ) \neq 0\}$ is a well-ordered subset of $\Gamma $,
that is, every nonempty subset of supp$(f)$ contains a minimal
element. For an arbitrary pair $f, g \in K _{0}((\Gamma ))$, the
sum $f + g$ and the product $f.g$ are defined by the rules $(f +
g)(\gamma ) = f(\gamma ) + g(\gamma )$, and $(f.g)(\gamma ) = \sum
_{\delta \in \Gamma } f(\gamma - \delta )g(\delta )$, for each
$\gamma \in \Gamma $ {\rm (}$f.g$ is correctly defined, since the
nonzero summands in the sum presenting $(f.g)(\gamma )$ are
finitely many, by known general properties of well-ordered subsets
of $\Gamma $, see \cite[Lemma~2.7.2]{E3}{\rm )}.}
\par
\medskip
The study of generalized power series fields as objects of
valuation theory, which form an important class of maximally
complete fields (containing $n$-fold iterated Laurent series
fields, see \cite[1.1.4]{TW}), dates back to Krull. His approach
to these fields is nowadays standard, and in this section, it is
used for characterizing the pairs $(K _{0}, \Gamma )$ for which $K
_{0}((\Gamma )) \in \Phi _{\rm Br}$, as follows:
\par
\medskip
\begin{prop}
\label{prop4.4} Assume that $K _{0}$ is a field and $\Gamma \neq 0$
is an ordered abelian group. Then the generalized power series field
$K = K _{0}((\Gamma ))$ is Brauer finite-dimensional if and only if
$K _{0}$ is virtually perfect and Brauer finite-dimensional, and the
quotient groups $\Gamma /p\Gamma $, $p \in \mathbb{P}$, are
finite.
\end{prop}
\par
\smallskip
\begin{proof} Let $v _{\Gamma }$ be the standard valuation of $K$
inducing on $K _{0}$ the trivial valuation. Then $(K, v _{\Gamma
})$ is maximally complete (cf. \cite[Theorem~18.4.1]{E3}) with
$v(K) = \Gamma $ and $\widehat K = K _{0}$ (see \cite[Sect. 2.8
and Example~4.2.1]{E3}), so Proposition \ref{prop4.4} can be
viewed as a special case of the following result.
\end{proof}
\par
\medskip
\begin{prop}
\label{prop4.5} Let $(K, v)$ be a Henselian field with {\rm
char}$(K) = {\rm char}(\widehat K)$, and in case {\rm
char}$(\widehat K) > 0$, suppose that $(K, v)$ is maximally
complete. Then $K \in \Phi _{\rm Br}$ if and only if $\widehat K$
is virtually perfect, and for each $p \in \mathbb{P}$, {\rm
abrd}$_{p}(\widehat K)$ and the group $v(K)/pv(K)$ are finite.
When this holds, $K$ is virtually perfect.
\end{prop}
\par
\medskip The latter assertion of Proposition \ref{prop4.5} is
implied by the former one and Lemma \ref{lemm4.1} if
char$(\widehat K) > 0$, and when char$(\widehat K) = 0$, it is
obvious. The former assertion follows from the next two lemmas and
the observations on absolute Brauer $p$-dimensions made at the end
of Section 3.
\par
\smallskip
\begin{lemm}
\label{lemm4.6} Let $(K, v)$ be a Henselian field and $p \neq {\rm
char}(\widehat K)$ a prime. Then {\rm abrd}$_{p}(K)$ is finite if
and only if so are {\rm abrd}$_{p}(\widehat K)$ and $v(K)/pv(K)$.
\end{lemm}
\par
\smallskip
\begin{proof}
We have abrd$_{p}(\widehat K) \le {\rm abrd}_{p}(K)$ (by
\cite[Theorem~2.8]{JW}, and \cite[Theorem~A.23]{TW}), so our
assertion can be deduced from \cite[Proposition~6.1, Theorem~5.9
and Remark~6.2]{Ch6} (or else, from \cite[(3.3) and
Theorem~2.3]{Ch6}).
\end{proof}
\par
\medskip
Lemma \ref{lemm4.6} characterizes those $\Phi _{\rm Br}$-fields
which are endowed with Henselian valuations whose residue fields
are of characteristic zero. Also, this lemma and the next one show
that a Henselian discrete valued field belongs to the class $\Phi
_{\rm Br}$ if and only if its residue field is a virtually perfect
$\Phi _{\rm Br}$-field.
\par
\medskip
\begin{lemm}
\label{lemm4.7} Let $(K, v)$ be a Henselian field with {\rm
char}$(\widehat K) = q > 0$. Then:
\par
{\rm (a)} $[\widehat K\colon \widehat K ^{q}]$ and $v(K)/qv(K)$
are finite, provided that so is {\rm Brd}$_{q}(K)$;
\par
{\rm (b)} The inequality {\rm abrd}$_{q}(K) < \infty $ holds,
provided that $[\widehat K\colon \widehat K ^{q}] < \infty $ and
$v$ is discrete or the following condition is satisfied:
\par
{\rm (i)} {\rm char}$(K) = q$ and $K$ is virtually perfect; this
occurs if {\rm char}$(K) = q$, $v(K)/qv(K)$ is finite and $(K, v)$
is maximally complete.
\end{lemm}
\par
\smallskip
\begin{proof}
Statement (a) is implied by \cite[Proposition~3.4]{Ch5}, so one
may assume that $[\widehat K\colon \widehat K ^{q}] = q ^{\mu }$
and $v(K)/qv(K)$ has order $q ^{\tau }$, for some integers $\mu
\ge 0$, $\tau \ge 0$. We prove Lemma \ref{lemm4.7} (b). Suppose
first that $v$ is discrete. Then the scalar extension map Br$(K)
\to {\rm Br}(K _{v})$ preserves Schur indices (cf.
\cite[Theorem~1]{Cohn}), which implies Brd$_{p}(K) \le {\rm
Brd}_{p}(K _{v})$ and abrd$_{p}(K) \le {\rm abrd}_{p}(K _{v})$,
for every $p \in \mathbb{P}$ (see also \cite[Ch. XII,
Corollary~4.5]{L}). This, applied to the case of $p = q$, shows
that it suffices to prove the inequality abrd$_{q}(K) < \infty $,
assuming that $K = K _{v}$. If char$(K) = 0$, then abrd$_{q}(K)
\le 2{\rm log}_{q}[\widehat K\colon \widehat K ^{q}] = 2\mu $, by
\cite[Theorem~2]{PS}; when char$(K) = q$, one proves that
abrd$_{q}(K) < \infty $ using the latter part of Lemma
\ref{lemm4.7} (b) (i) and the fact that $(K, v)$ is maximally
complete.
\par
In the rest of our proof, we drop the assumption that $v$ is
discrete and suppose that char$(K) = q$. Then it follows from
Lemma \ref{lemm3.7} (b) that abrd$_{q}(K) < \infty $ whenever
$[K\colon K ^{q}] < \infty $. Note finally that if $(K, v)$ is
maximally complete, then $[K\colon K ^{q}] = q ^{\mu + \tau }$.
Indeed, Lemma \ref{lemm4.1} (c), applied to the valued subfield
$(K ^{q}, v _{q})$ of $(K, v)$, implies $q ^{\mu + \tau }$ is an
upper bound of the degrees of finite extensions of $K ^{q}$ in
$K$. This yields $[K\colon K ^{q}] \le q ^{\mu + \tau }$, which
allows to conclude that $[K\colon K ^{q}] = q ^{\mu + \tau }$ and
so completes the proof of Lemma \ref{lemm4.7} (b).
\end{proof}
\par
\medskip
\section{\bf Proof of Theorem \ref{theo2.3}}
\par
\smallskip The former part of Theorem \ref{theo2.3} (a) is a
special case of Theorem \ref{theo3.1}, since the class of $C
_{m}$-fields is closed under the formation of algebraic
extensions, for each $m$, and by Matzri's theorem, every $C
_{m}$-field $E$ is virtually perfect with Brd$_{p}(E) \le p ^{m-1}
- 1$, for all $p \in \mathbb{P}_{E}$. Similarly, the latter part
of Theorem \ref{theo2.3} (a) is a special case of the former one,
since ddim$(E _{0}) = 1$ if $E _{0}$ is a finite field, ddim$(E
_{0}) \le 2$ if $E _{0}$ is a perfect PAC-field, and by
Lang-Nagata-Tsen's theorem, the class of fields of finite
Diophantine dimensions is closed under taking field extensions of
finite transcendence degree.
\par
\smallskip
We turn to the proof of Theorem \ref{theo2.3} (b). The general
part of this result is contained in the following characterization
of $m$-local $\Phi _{\rm Br}$-fields.
\par
\medskip
\begin{prop}
\label{prop5.1} Let $K _{m}$ be an $m$-local field with an $m$-th
residue field $K _{0}$, for some $m \in \mathbb{N}$. Then $K _{m}$
is a Brauer finite-dimensional field if and only if $K _{0}$ is
virtually perfect with {\rm abrd}$_{p}(K _{0}) < \infty $, for
every $p \in \mathbb{P}_{K_{0}}$.
\end{prop}
\par
\smallskip
\begin{proof}
Arguing by induction on $m$, and using Lemma \ref{lemm4.6}, one
obtains that, for any $p \in \mathbb{P}_{K_{0}}$, abrd$_{p}(K
_{m}) < \infty \leftrightarrow {\rm abrd}_{p}(K _{0}) < \infty $.
Thus our proof is complete in the case of char$(K _{0}) = 0$, so
we assume that char$(K _{0}) = q > 0$. Let $K _{m-j}$ be the
$j$-th residue field of $K _{m}$, for $j = 1, \dots , m$. Then
Brd$_{q}(K _{m-j+1}) \ge {\rm Brd}_{q}(K _{m-j})$ (see
\cite[Theorem~2.8]{JW}) which implies in
\par\vskip0.04truecm\noindent
conjunction with Lemma \ref{lemm4.3} that abrd$_{q}(K _{m-j+1}) \ge
{\rm abrd}_{q}(K _{m-j})$, for each $j$.
Since, by Lemma \ref{lemm4.7} (a), Brd$_{q}(K _{1}) = \infty $ if $[K
_{0}\colon K _{0} ^{q}] = \infty $, this allows to assume in the rest
of the proof of Proposition \ref{prop5.1} that $[K _{0}\colon K _{0}
^{q}] < \infty $. Suppose first that char$(K _{m}) = q$. Then, by the
Hasse-Schmidt theorem (see \cite[Proposition~A.5.3]{GiSz}, or
\cite[page~110]{E3}, and further references there), $K _{m-j+1}$ is
isomorphic to the Laurent series field $K _{m-j}((X _{m-j+1}))$, $j =
1, \dots , m$, which implies the degree $[K _{m-j+1}\colon K _{m-j+1}
^{q}]$ equals $q ^{m-j+1}.[K _{0}\colon K _{0} ^{q}]$, for each $j$.
Now it follows from Lemma \ref{lemm3.7} (b) that Brd$_{q}(K _{m}
^{\prime }) \le {\rm log}_{q}[K _{m}\colon K _{m} ^{q}]$, for every
finite extension $K _{m} ^{\prime }/K _{m}$, which completes the
proof of Proposition \ref{prop5.1} in the case of char$(K _{m}) = q$.
It remains to be seen that abrd$_{q}(K _{m}) < \infty $ when char$(K
_{m}) = 0$. Then char$(K _{m-\mu + 1}) = 0$ and char$(K _{m-\mu }) =
q$, for some $\mu \le m$, and by the previous part of our proof, $[K
_{m-\mu }\colon K _{m-\mu } ^{q}] < \infty $; hence, by
\cite[Theorem~2]{PS}, abrd$_{q}(K _{m-\mu + 1}) < \infty $, which
enables one to deduce from Lemma \ref{lemm4.6} that abrd$_{q}(K _{m})
< \infty $, as required.
\end{proof}
\par
\smallskip
\begin{coro}
\label{coro5.2} The $m$-fold iterated Laurent series field $K
_{0}((Z _{1})) \dots ((Z _{m}))$, where $m \in \mathbb{N}$, is
Brauer finite-dimensional if and only if the coefficient field $K
_{0}$ is virtually perfect and Brauer finite-dimensional.
\end{coro}
\par
\smallskip
\begin{proof}
This is a special case of Propositions \ref{prop4.4} and
\ref{prop5.1}.
\end{proof}
\par
\medskip
To prove the latter part of Theorem \ref{theo2.3} (b) we show that
a field $K _{0}$ is Brauer finite-dimensional, provided that
ddim$(K _{0}) < \infty $ or $K _{0}$ is a finitely-generated
extension of a PAC field $E$. When ddim$(K _{0}) < \infty $ or
char$(E) = 0$, this has already been done, so we assume that $K
_{0}$ is a finitely-generated extension of some PAC field $E$ with
char$(E) = q > 0$; in the setting of Theorem \ref{theo2.3}, this
requires that $E$ be virtually perfect. Identifying $E _{\rm ins}$
with its $E$-isomorphic copy in an algebraic closure $\overline K
_{0}$ of $K _{0}$, and observing that $E _{\rm ins}$ is PAC (cf.
\cite{FJ}), one concludes that $E _{\rm ins}K _{0}/E _{\rm ins}$
is a finitely-generated field extension, ddim$(E _{\rm ins})$ and
ddim$(E _{\rm ins}K _{0})$ are finite, $E _{\rm ins}K _{0} \in I(K
_{0,{\rm ins}}/K _{0})$, and $(E _{\rm ins}K _{0}) _{\rm ins} = K
_{0,{\rm ins}}$. Therefore, Lemma \ref{lemm3.7} yields Brd$_{p}(K
_{0}) = {\rm Brd}_{p}(E _{\rm ins}K _{0}) = {\rm Brd}_{p}(K
_{0,{\rm ins}}) < \infty $, for all $p \in \mathbb{P} _{E}$. Since
$K _{0}$ is virtually perfect, these equalities show that it is
Brauer finite-dimensional and so complete the proof of Theorem
\ref{theo2.3} (b).
\par
\medskip
Applying Matzri's theorem and using repeatedly Greenberg's theorem
(cf. \cite{Gr}), one illustrates the former part of Theorem
\ref{theo2.3} (a) as follows:
\par
\medskip
\begin{prop}
\label{prop5.3} Let $K _{0}$ be a virtually perfect field, $K
_{m}$ the $m$-fold Laurent series field $K _{0}((Z _{1})) \dots
((Z _{m}))$, for some $m \in \mathbb{N}$, and $K/K _{m}$ a field
extension. Suppose that {\rm ddim}$(K _{0})$ and {\rm trd}$(K/K
_{m})$ are finite. Then {\rm ddim}$(K)$ is finite.
\end{prop}
\par
\smallskip
The following lemma is used for proving Theorem \ref{theo2.3} (c)
in the case where $K _{0}$ is a virtually perfect PAC field.
\par
\smallskip
\begin{lemm}
\label{lemm5.4} Let $F$ be a finitely-generated extension of a
perfect {\rm PAC} field $E$ with {\rm trd}$(F/E) = 1$. Then {\rm
Brd}$_{p}(F) = {\rm abrd}_{p}(F) = {\rm cd}_{p}(\mathcal{G}_{E})
\le 1\colon p \in \mathbb{P}$.
\end{lemm}
\par
\smallskip
\begin{proof}
Here we identify $E _{\rm sep}$ with its $E$-isomorphic copy in an
algebraic closure $\overline F$ of $F _{\rm sep}$. It is known
that cd$_{p}(\mathcal{G}_{E}) \le 1$, for all $p \in \mathbb{P}$
(cf. \cite[Corollary~11.6.8]{FJ}); also, we have Brd$_{p}(F) \ge
1$ whenever cd$_{p}(\mathcal{G}_{E}) > 0$ (see
\cite[Proposition~5.9]{Ch4}, and the reference to its proof). On
the other hand, it follows from Galois theory (and Galois
cohomology, see \cite[Ch. I, 3.3]{S1}) that
cd$_{p}(\mathcal{G}_{E}) = 0$ if and only if $p \nmid [E
_{1}\colon E]$, for any finite extension $E _{1}$ of $E$ in $E
_{\rm sep}$. Observing that $E _{\rm sep}F ^{\prime }/F ^{\prime
}$ is a Galois extension with $\mathcal{G}(E _{\rm sep}F ^{\prime
}/F ^{\prime })$ isomorphic to an open subgroup of
$\mathcal{G}_{E}$, for each finite extension $F ^{\prime }$ of $F$
in $\overline F$ (cf. \cite[Ch. VI, Theorem~1.12]{L}), one proves
that the equality cd$_{p}(\mathcal{G}_{E}) = 0$ implies $p \nmid
[F _{1} ^{\prime }\colon F ^{\prime }]$, for any finite extension
$F _{1} ^{\prime }$ of $F ^{\prime }$ in $E _{\rm sep}F ^{\prime
}$. Therefore, the scalar extension map of $s(F ^{\prime })$ into
$s(E _{\rm sep}F ^{\prime })$ induces an embedding of Br$(F
^{\prime }) _{p}$ into Br$(E _{\rm sep}F ^{\prime }) _{p}$ (apply,
e.g., \cite[Sect.~13.4, Proposition~(vi)]{P}). Since, by Tsen's
theorem (cf. \cite[Sect.~19.4]{P}), ddim$(E _{\rm sep}F ^{\prime
}) \le 1$, this yields Br$(F ^{\prime }) _{p} = \{0\}$, proving
that abrd$_{p}(F) = 0$. For the proof of Lemma \ref{lemm5.4}, it
remains to be shown that if the set $\mathbb{P}_{1} = \{p \in
\mathbb{P}\colon {\rm cd}_{p}(\mathcal{G}_{E}) = 1\}$ is nonempty,
then abrd$_{p}(F) = 1$, for all $p \in \mathbb{P}_{1}$. By
Albert's criterion (cf. \cite[Ch. XI, Theorem~3]{A1}), it suffices
to prove that if $Y/F$ is a finite extension, then deg$(D _{p}) =
p$, for each $D _{p} \in d(Y)$ with exp$(D _{p}) = p$. Since
algebraic extensions of $E$ are perfect and PAC fields (cf.
\cite[Corollary~11.2.5]{FJ}), one may assume without loss of
generality that $E$ has no proper algebraic extension in $Y$. Fix
a full system of representatives $V(Y)$ of the set of equivalence
classes of nontrivial valuations of $Y$, that are trivial on $E$.
It is known that abrd$_{p}(E) = 0\colon p \in \mathbb{P}$, and
every $v \in V(Y)$ is discrete with a residue field depending on
$v$ that is a finite extension of $Y$ (cf. \cite[Theorem~11.6.4
and Example~2.2.1~(b)]{FJ}, and \cite[Ch. XII, Corollary~4.5]{L}).
Hence, by Witt's theorem and Lemma \ref{lemm4.2} (c), every $D
_{v} \in d(Y _{v})$ is a cyclic NSR-algebra (see also
\cite[Proposition~8.75]{TW}, and \cite[page~133 and Lemmas~5.1,
5.14]{JW}). This ensures that abrd$_{p}(Y _{v}) \le 1\colon p \in
\mathbb{P}$, for all $v \in v(Y)$ (e.g., by
\cite[Theorem~5.15~(a)]{JW} or \cite[Propositions~5.3~(b),
5.4]{Ch6}). Observe finally that the Hasse Principle for algebras
from $s(Y)$, and Cohn's theorem (cf. \cite[Theorem~3.4]{Ef} and
\cite[Theorem~1]{Cohn}, respectively) imply the scalar extension
maps $s(Y) \to s(Y _{v})$, $v \in v(Y)$, give rise to an embedding
of Br$(Y)$ into the direct sum $\oplus _{v \in v(Y)} {\rm Br}(Y
_{v})$. Therefore, one obtains from Grunwald-Wang's theorem (in
fact, from its generalization in \cite{LR}), by the method of
proving \cite[(32.19)]{Re}, that $D _{p}$ has a splitting field $Y
_{p} ^{\prime }$ which is a degree $p$ cyclic extension of $Y$.
Thus it turns out that $Y _{p} ^{\prime }$ is $E$-isomorphic to a
maximal subfield of $D _{p}$ and deg$(D _{p}) = p$ (cf.
\cite[Sect.~13.4, Corollary]{P}).
\end{proof}
\par
\medskip
{\it Proof of Theorem \ref{theo2.3} {\rm (c)}}. It follows from
Hasse-Schmidt's theorem that the conditions of Proposition
\ref{prop5.3} are equivalent to the one that $K$ is an extension
of an $m$-local field $K _{m}$ with a virtually perfect $m$-th
residue field $K _{0}$, such that ddim$(K _{0}) < \infty $,
trd$(K/K _{m}) < \infty $, and char$(K _{m}) = {\rm char}(K
_{0})$. In other words, the assertion of Theorem \ref{theo2.3} (c)
in the special case of char$(K _{m}) = {\rm char}(K _{0})$ is
generalized by Proposition \ref{prop5.3}, so one may assume
further that char$(K _{m}) = 0$ and char$(K _{0}) = q > 0$. In
view of Theorem \ref{theo3.1}, Proposition \ref{prop5.1} and
\cite[(1.3)]{Ch3}, it suffices to show that $K \in \Phi _{\rm Br}$
when trd$(K/K _{m}) = 1$ and $K/K _{m}$ is finitely-generated. Let
$K _{0} ^{\prime }/K _{0}$ be a finitely-generated extension with
trd$(K _{0} ^{\prime }/K _{0}) = 1$. Then it follows from Lemmas
\ref{lemm3.7} and \ref{lemm5.4} that if the field $K _{0}$ is
virtually perfect and PAC, then $K _{0} ^{\prime }$ is a virtually
perfect $\Phi _{\rm Br}$-field with Brd$_{p}(K _{0} ^{\prime })
\le 1$, for all $p \in \mathbb{P}_{K _{0}}$. When $K _{0}$ is a
finitely-generated extension of a finite or algebraically closed
field $E$, such that trd$(K _{0}/E) \le 1$, we have Brd$_{p}(K
_{0} ^{\prime }) = e \le 2\colon p \in \mathbb{P}_{K _{0}}$, where
$e$ is determined as follows: $e = 0$ if $K _{0} = E$ is
algebraically closed (by Tsen's theorem); $e = 2$ if $E$ is finite
and trd$(K _{0}/E) = 1$ (see \cite{Li15}); $e = 1$ if $E$ is
finite and trd$(K _{0}/E) = 0$ (cf. \cite[(32.19)]{Re}); $e = 1$
if $E$ is algebraically closed and trd$(K _{0}/E) = 1$ (cf.
\cite{Jong} and \cite{Li08}). Therefore, the assertion that $K \in
\Phi _{\rm Br}$ can be deduced from \cite[Theorem~3]{PS}, applied
as in the proof of Proposition \ref{prop5.1}, and the
Harbater-Hartmann-Krashen-Lieblich theorem (see
\cite[Corollaries~5.7, 5.8]{HHKr} and \cite[Corollary~1.3]{Li11}).
Theorem \ref{theo2.3} is proved.
\par
\medskip
It has been proved in \cite{Ch7} that the class $\Phi _{\rm Br}$
contains the function field $K _{C}$ of an arbitrary algebraic
curve $C$ defined over an $m$-local field $K _{m}$ with an $m$-th
residue field $K _{0}$ of finite Diophantine dimension; hence, by
Theorem \ref{theo3.1}, central division LFD-algebras over $K _{C}$
are NLF. The result has been obtained as a consequence of the
Lang-Nagata-Tsen theorem, Matzri's theorem and a version of the
Harbater-Hartmann-Krashen-Lieblich theorem, proved in \cite{Ch7}.
This version facilitates the verification of whether a given field
lies in $\Phi _{\rm Br}$, and thereby, the application of Theorem
\ref{theo3.1}.
\par
\medskip
\begin{rema}
\label{rema5.5} A well-known conjecture, originally proposed by M.
Artin for $d = 2$, predicts that if $K$ is a $C _{d}$-field, for
some $d > 0$, then the dimensions {\rm Brd}$_{p}(K)\colon p \in
\mathbb{P}$, have an upper bound $d - 1$ (see
\cite[Sect.~4]{ABGV}). If the conjecture holds, it will become
possible, as in the proof of Theorem \ref{theo2.3} {\rm (c)}, to
deduce the above-noted result of \cite{Ch7} directly from the
Harbater-Hartmann-Krashen-Lieblich theorem. Note that Lemma
\ref{lemm5.4} agrees with Artin's original conjecture, especially,
if char$(E) = 0$ or $E$ is perfect with $\mathcal{G}_{E}$ a
pro-$\ell $ group, for some $\ell \in \mathbb{P}$, since then
ddim$(E) \le 1$ (cf. \cite{Koll} and \cite[Sect. 3]{Witt}), and
therefore, ddim$(F) \le 2$.
\end{rema}
\par
\medskip
Theorems \ref{theo3.1} and \ref{theo2.3}, as well as Matzri's
theorem and the Lang-Nagata-Tsen theorem, attract interest in the
open question of whether the class of virtually perfect $\Phi
_{\rm Br}$-fields is closed under taking finitely-generated field
extensions. Using Galois theory, Sylow's theorems and
\cite[Lemma~3.5]{Ch2}, one concludes that the answer to this
question will be affirmative if it turns out that
ddim$(\mathcal{E} _{p}) < \infty $, for each $p \in \mathbb{P}$
and every field $\mathcal{E} _{p}$ with char$(\mathcal{E} _{p})
\neq p$, $\mathcal{E} _{p}(p) = \mathcal{E} _{p,{\rm sep}}$ and
abrd$_{p}(\mathcal{E} _{p}) < \infty $. The answer is presently
unknown in the case where $\mathcal{E} _{p}$ is the fixed field of
a Sylow pro-$p$ subgroup of $\mathcal{G}_{\mathbb{Q}}$ or
$\mathcal{G}_{\mathbb{Q}_{\ell }}$, for some $\ell \in
\mathbb{P}$.
\par
\smallskip
\section{\bf Preparation for the proof of Proposition \ref{prop2.4}}
\par
\medskip
In this Section we collect results on finite-dimensional central
division algebras over Henselian fields, which are used in Section
7 for proving Proposition \ref{prop2.4}. We consider mainly tame
division algebras central over a strictly Henselian field. Our
starting point is the following lemma.
\par
\medskip
\begin{lemm}
\label{lemm6.1} Let $(K, v)$ be a Henselian field, $D \in d(K)$ a
division $K$-algebra of degree not divisible by {\rm
char}$(\widehat K)$, and $\Delta _{0}$ an inertial central
$K$-subalgebra of $D$. Then $\Delta _{0}$ is included in an
inertial $K$-subalgebra $\Delta _{1}$ of $D$ with $\widehat \Delta
_{1} = \widehat D$.
\end{lemm}
\par
\smallskip
\begin{proof}
It is known that $D = \Delta _{0} \otimes _{K} C _{0}$, where $C
_{0} = C _{D}(\Delta _{0})$, and also, that $C _{0}$ contains as a
$K$-subalgebra an inertial lift $I _{1}$ of $\widehat C _{0}$ over
$K$, which is uniquely determined, up-to $K$-isomorphism (cf.
\cite[Theorem~4.4.2]{He}, and \cite[Theorem~2.9]{JW},
respectively). Clearly, $\Delta _{0}$ is a subalgebra of the
$K$-algebra $\Delta _{1} = \Delta _{0} \otimes _{K} I _{1}$, and
it is not difficult to see, using \cite[Theorems~3.1 and 2.8]{JW},
that if $Z _{1} = Z(I _{1})$, then $\Delta _{0} \otimes _{K} Z _{1}$
and $\Delta _{1} \cong (\Delta _{0} \otimes _{K} Z _{1}) \otimes _{Z
_{1}} I _{1}$ are inertial division $K$-algebras central over $Z
_{1}$. Moreover, it follows that $Z(\widehat C _{0}) = \widehat Z
_{1}$
\par\vskip0.04truecm\noindent
and $[\widehat \Delta _{1}] = [(\widehat \Delta _{0} \otimes
_{\widehat K} \widehat Z _{1}) \otimes _{\widehat Z _{1}} \widehat
C _{0}] = [\widehat \Delta _{0} \otimes _{\widehat K} \widehat C
_{0}]$ in Br$(\widehat Z _{1})$. Since
\par\vskip0.04truecm\noindent
char$(\widehat K) \nmid {\rm deg}(D)$, one finally obtains from
\cite[Corollary~6.8]{JW}, that
\par\vskip0.04truecm\noindent
$v(C _{0}) = v(D)$, $\widehat Z _{1} = Z(\widehat D)$, and
$[\widehat D] = [\widehat \Delta _{0} \otimes _{\widehat K}
\widehat C _{0}] = [\widehat \Delta _{1}]$, which proves Lemma
\ref{lemm6.1}.
\end{proof}
\par
\smallskip
Under the hypothesis that $(K, v)$ is a strictly Henselian field,
and $D _{1}, \dots , D _{s}$ are $K$-algebras with $D _{i} \in
d(K)$ and char$(\widehat K) \nmid [D _{i}\colon K]$, $i = 1, \dots
, s$, for some integer $s \ge 2$, it follows from Lemma
\ref{lemm4.2} (a), Draxl's lemma and Morandi's
\par\vskip0.04truecm\noindent
theorem (see \cite[Lemma~3]{Dr1} and \cite[Theorem~1]{Mo}) that
\par\vskip0.04truecm\noindent
$D _{1} \otimes _{K} \dots \otimes _{K} D _{s} \in d(K)$
if and only if the sum of the subgroups
\par\vskip0.04truecm\noindent
$v(D _{1})/v(K), \dots , v(D _{s})/v(K)$ of $\overline {v(K)}/v(K)$
is direct. This plays a major role in the proof of the following lemma.
\par
\medskip
\begin{lemm}
\label{lemm6.2} Let $(K, v)$ be a strictly Henselian field, and $D
_{3} \in d(K)$ be an algebra isomorphic to $D _{1} \otimes _{K} D
_{2}$, for some $D _{i} \in d(K)$, $i = 1, 2$. Suppose that {\rm
char}$(\widehat K) \nmid {\rm deg}(D)$, fix a divisor $\nu \in
\mathbb{N}$ of {\rm exp}$(D _{1})$, {\rm exp}$(D _{2})$ and {\rm
exp}$(D)$, and take algebras $\Delta _{j} \in d(K)$, $j = 1, 2,
3$, so that $[\Delta _{j}] = \nu [D _{j}]$, for each index $j$.
Then $\Delta $ and $\Delta _{1} \otimes _{K} \Delta _{2}$ are
isomorphic $K$-algebras.
\end{lemm}
\par
\medskip
\begin{proof}
The condition that $D _{1} \otimes _{K} D _{2} \in d(K)$, the
assumptions on $(K, v)$ and deg$(D)$, and Draxl's lemma ensure
that $v(D _{1}) \cap v(D _{2}) = v(K)$ and the group $v(D)/v(K)$
is isomorphic to the direct sum $v(D _{1})/v(K) \oplus \ v(D
_{2})/v(K)$. Applying \cite[Proposition~6.9]{JW}, one obtains
further that $v(\Delta _{i}) = \nu v(D _{i}) + v(K)$, $i = 1, 2$,
which yields $v(\Delta _{1}) \cap v(\Delta _{2}) = v(K)$. As
$\widehat \Delta _{1} = \widehat \Delta _{2} = \widehat K$ and
char$(\widehat K)$ does not divide deg$(\Delta _{1}){\rm
deg}(\Delta _{2})$, this combined with \cite[Theorem~1]{Mo},
proves that $\Delta _{1} \otimes _{K} \Delta _{2} \in d(K)$. Now
our assertion follows from the condition that $\Delta \in d(K)$
and the equality $[\Delta _{1} \otimes _{K} \Delta _{2}] = [\Delta
]$, in Br$(K)$.
\end{proof}
\par
\medskip
It is known that if $(K, v)$ is a Henselian field, then every tame
and totally ramified extension $T/K$ is radical, i.e. $T$ is
generated over $K$ by roots of elements of $K ^{\ast }$ of degrees
determined by the group $v(T)/v(K)$ (see
\cite[Proposition~A22]{TW}). When $n \ge 2$ is an integer and
char$(\widehat K) \nmid n$, our next lemma gives a necessary and
sufficient condition for an extension $K _{n}$ of $K$ obtained by
adjunction of an $n$-th root of an element $\alpha \in K ^{\ast }$
to be totally ramified of degree $n$. It shows that $K _{n}$ is
uniquely determined, up-to a $K$-isomorphism, by the group $(K
_{n} ^{\ast n} \cap K ^{\ast })/K ^{\ast n}$.
\par
\medskip
\begin{lemm}
\label{lemm6.3} Let $(K, v)$ be a Henselian field, $n \ge 2$ an
integer not divisible by {\rm char}$(\widehat K)$, $K _{n}$ an
extension of $K$ in $K _{\rm sep}$ generated by an $n$-th root
$\alpha _{n}$ of an element $\alpha \in K ^{\ast }$, and $N(K
_{n}/K)$ the norm group of $K _{n}/K$. Then $K _{n}/K$ is totally
ramified of degree $n$ if and only if the coset $v(\alpha ) +
nv(K)$ has order $n$ as an element of the group $v(K)/nv(K)$. Such
being the case, the following holds:
\par\vskip0.04truecm\noindent
$\ $ {\rm (a)} $K _{n} ^{\ast }$, $K _{n} ^{\ast n} \cap K ^{\ast
}$ and $N(K _{n}/K)$ coincide with the sets $\cup _{u=0} ^{n-1}
\alpha _{n} ^{u}(K ^{\ast }.K _{n} ^{\ast n})$,
\par\vskip0.04truecm\noindent
$\cup _{u=0} ^{n-1} \alpha ^{u}K ^{\ast n}$ and $\cup _{u=0}
^{n-1} ((-1) ^{(n-1)}.\alpha ) ^{u}K ^{\ast n}$, respectively; in
particular, the group
\par\vskip0.04truecm\noindent
$(K _{n} ^{\ast n} \cap K ^{\ast })/K ^{\ast n}$ is cyclic of
order $n$, and is generated by the coset $\alpha K ^{\ast n}$.
\par\noindent
$\ $ {\rm (b)} $N(K _{n}/K) = K _{n} ^{\ast n} \cap K ^{\ast }$ if
$2 \nmid n$ or $\widehat K$ contains a primitive $2n$-th root of
$1$.
\end{lemm}
\par
\smallskip
\begin{proof}
Our assumptions show that $[K _{n}\colon K] \le n$, and equality
holds if and only if the polynomial $g _{\alpha }(X) = X ^{n} -
\alpha \in K[X]$ is irreducible over $K$. Note further that the
coset $v(\alpha ) + nv(K)$ has order $n$ in $v(K)/nv(K)$ if and
only if the coset $(1/n)v(\alpha ) + v(K)$ has order $n$ in
$\overline {v(K)}/v(K)$ (cf. \cite[Lemma~7.70]{TW}). As $v _{K
_{n}}(\alpha _{n}) = (1/n)v(\alpha )$, the fulfillment of this
condition requires $e(K _{n}/K) \ge n$, which enables one to deduce
from Ostrowski's theorem and the inequality $[K _{n}\colon K] \le n$
that $[K _{n}\colon K] = e(K _{n}/K) = n$ and
\par\vskip0.04truecm\noindent
$v(K _{n}) = \cup _{i=0} ^{n-1} (v _{K _{n}}(\alpha _{n} ^{i}) +
v(K))$. Let now $[K _{n}\colon K] = n$ and the order of
\par\vskip0.04truecm\noindent
$v(\alpha ) + nv(K)$ in $v(K)/nv(K)$ be less than $n$. Then it
follows that
\par\vskip0.04truecm\noindent
$v(\alpha ) \in pv(K)$, for some $p \in \mathbb{P}$ dividing $n$.
This implies $\alpha = \pi ^{p}\alpha _{0}$, for some $\pi , \alpha
_{0} \in K ^{\ast }$, chosen so that $v(\alpha _{0}) = 0$. Observing
that $\alpha _{n} ^{n} = \alpha $ and $g _{\alpha }(X)$ is
irreducible over $K$, one concludes that $\alpha _{0} \in (K _{n}
^{\ast p} \setminus K ^{\ast p})$. As $(K, v)$ and $(K _{n}, v _{K
_{n}})$ are Henselian, (4.2) (a) and this fact indicate that the
residue $\hat \alpha _{0}$ lies in the set $\widehat K _{n} ^{\ast p}
\setminus \widehat K ^{\ast p}$; hence, $\widehat K _{n} \neq
\widehat K$ and $e(K _{n}/K) < [K _{n}\colon K]$, which completes the
proof of our assertion about $K _{n}/K$.
\par
In the rest of our proof, one may assume that $e(K _{n}/K) = [K
_{n}\colon K] = n$. Then $\widehat K _{n} = \widehat K$, so it
follows from the above description of $v(K _{n})$ that each
$\lambda _{n} \in K _{n} ^{\ast }$ is presentable as a product
$\lambda _{n} = \alpha _{n} ^{u}\lambda \lambda _{0}$, for some
integer $u = u(\lambda _{n})$ and elements $\lambda \in K ^{\ast
}$, $\lambda _{0} \in \nabla _{0}(K _{n})$ with $0 \le u < n$.
Therefore, by (4.2) (a) (and the condition that char$(\widehat K)
\nmid n$), $\lambda _{0} \in K _{n} ^{\ast n}$, which implies the
image of
\par\vskip0.04truecm\noindent
$\lambda _{n}$ under the norm map $N _{K _{n}/K}\colon K _{n}
^{\ast } \to K ^{\ast }$ equals $(-1) ^{(n-1)u}\alpha ^{u}\lambda
^{n}\tilde \lambda _{0} ^{n}$, where $\tilde \lambda _{0} ^{n} = N
_{K _{n}/K}(\lambda _{0})$. These calculations prove the
statements of Lemma \ref{lemm6.3}~(a) about $K _{n} ^{\ast }$ and
$N(K _{n}/K)$. Since $\widehat K _{n} = \widehat K$ and $\nabla
_{0}(K _{n}) \cap K ^{\ast } \subseteq \nabla _{0}(K) \subseteq K
^{\ast n}$,
\par\vskip0.04truecm\noindent
they also show that $K _{n} ^{\ast n} \cap K ^{\ast } = \cup
_{u=0} ^{n-1} \alpha ^{u}K ^{\ast n}$, and the coset $\alpha K
^{\ast }$, viewed as
\par\vskip0.04truecm\noindent
an element of $K ^{\ast }/K ^{\ast n}$, generates the quotient $(K
_{n} ^{\ast n} \cap K ^{\ast })/K ^{\ast n}$. Now the latter part
of Lemma \ref{lemm6.3} (a) becomes obvious. As to Lemma
\ref{lemm6.3}~(b), it follows from the former part of Lemma
\ref{lemm6.3}~(a), since $(-1) ^{(n-1)} = 1$ if $2 \nmid n$, and
$-1 \in K ^{\ast n}$ if $2 \mid n$ and $\widehat K$ contains a
primitive $2n$-th root of unity.
\end{proof}
\par
\medskip
The following lemma is a version of \cite[Theorem~1 and
Lemma~3]{Dr1}, which describes the structure of central division
tame algebras $T$ over a strictly Henselian field, as well as its
relations with $v(T)$. The former part of the lemma characterizes
the division algebras among the symbol algebras over the
considered field. The lemma is known (cf.
\cite[Proposition~2]{PY}, \cite[Corollary~2.6]{JW1},
\cite[Corollary~6.10]{JW}, \cite[Sect. 7.4]{TW}, and for an analog
about graded algebras, see \cite[Proposition~7.57]{TW}). For
convenience of the reader, we prove here its former part in the
traditional framework of valuation theory and radical field
extensions.
\par
\medskip
\begin{lemm}
\label{lemm6.4} Let $(K, v)$ be a strictly Henselian field, $n \ge
2$ an integer not divisible by {\rm char}$(\widehat K)$, and $u
_{n}$ a primitive $n$-th root of unity in $K$. Then:
\par
{\rm (a)} The class $d(K)$ contains the symbol $K$-algebra $\Delta
_{n} = K(\alpha , \beta ; u _{n}) _{n}$, where $\alpha , \beta \in
K ^{\ast }$, if and only if the cosets $v(\alpha ) + nv(K)$ and
$v(\beta ) + nv(K)$ generate a subgroup of $v(K)/nv(K)$ of order
$n ^{2}$; if $\Delta _{n} \in d(K)$, then {\rm exp}$(\Delta _{n})
= n$, the quotient group $v(\Delta _{n})/v(K)$ is generated by the
cosets $(1/n)v(\alpha ) + v(K)$ and $(1/n)v(\beta ) + v(K)$, and
it is isomorphic to
\par\vskip0.032truecm\noindent
$(\mathbb{Z}/n\mathbb{Z}) ^{2} = (\mathbb{Z}/n\mathbb{Z}) \oplus
(\mathbb{Z}/n\mathbb{Z})$.
\par
{\rm (b)} An algebra $D \in d(K)$ of degree prime to {\rm
char}$(\widehat K)$ is either a symbol one or deg$(D) > {\rm
exp}(D)$ and $D$ is isomorphic to the $K$-algebra $\otimes _{i=1}
^{\nu } S _{i}$, where $\otimes = \otimes _{K}$, $\nu \ge 2$ is an
integer, and $S _{i} = K(\alpha _{i}, \beta _{i}; u _{n} ^{(n/n
_{i})}) _{n _{i}}$, for some $\alpha _{i}, \beta _{i} \in K ^{\ast
}$ and each $i$; in the latter case, deg$(D) = n = \prod _{i=1}
^{\nu } n _{i}$, {\rm exp}$(D)$ is the least common multiple of $n
_{1}, \dots , n _{\nu }$, and there are group isomorphisms
\par\vskip0.032truecm\noindent
$v(D)/v(K) \cong G \oplus G$, $G \cong \oplus _{i=1} ^{\nu }
(\mathbb{Z}/n _{i}\mathbb{Z})$, and $v(S _{i})/v(K) \cong
(\mathbb{Z}/n _{i}\mathbb{Z}) ^{2}$, for $i = 1, \dots , \nu $.
\end{lemm}
\par
\smallskip
\begin{proof}
Lemma \ref{lemm6.4} (b) follows from Lemma \ref{lemm6.4} (a) and
\cite[Theorem~1 and Lemma~3]{Dr1}, so we prove only Lemma
\ref{lemm6.4} (a). Since, by definition, $\Delta _{n}$ is
generated by elements $x _{n}$ and $y _{n}$ subject to the
relations $y _{n}x _{n} = u _{n}x _{n}y _{n}$, $x _{n} ^{n} =
\alpha $, and $y _{n} ^{n} = \beta $, it is easily verified that
$\Delta _{n} \in s(K)$, deg$(\Delta _{n}) = n$, and the system $x
_{n} ^{i}y _{n} ^{j}\colon 0 \le i, j < n$, is a basis of $\Delta
_{n}$. This implies the rings $K[x _{n}]$ and $K[y _{n}]$ are
$K$-algebras isomorphic to the polynomial quotient rings $K[X]/(X
^{n} - \alpha )$ and $K[X]/(X ^{n} - \beta )$, respectively; in
particular, $[K[x _{n}]\colon K] = [K[y _{n}]\colon K] = n$.
Moreover, it follows that $K[x _{n}]$ and $K[y _{n}]$ are fields
unless $K[x _{n}]$ or $K[y _{n}]$ has nontrivial zero-divisors,
and both rings are fields if and only if the polynomials $X ^{n} -
\alpha $ and $X ^{n} - \beta $ are irreducible over $K$.
Therefore, it suffices to prove Lemma \ref{lemm6.4}~(a) under the
assumption that $X ^{n} - \alpha $ and $X ^{n} - \beta $ are
irreducible over $K$. Since, by (4.2) (b), $K$ contains a
primitive $n'$-th root of unity whenever $n' \in \mathbb{N}$ and
char$(\widehat K) \nmid n'$, this means that $\alpha \notin K
^{\ast p}$ and $\beta \notin K ^{\ast p}$, for any $p \in
\mathbb{P}$ dividing $n$ (cf. \cite[Theorem~9.1]{L}), so it turns
out that the cosets $\alpha K ^{\ast n}$ and $\beta K ^{\ast n}$
have order $n$ in $K ^{\ast }/K ^{\ast n}$. Hence, by Kummer
theory, $K[x _{n}]$ and $K[y _{n}]$ are cyclic extensions of $K$,
which are totally ramified, since $\widehat K _{\rm sep} =
\widehat K$ and char$(\widehat K) \nmid n$. Also, by Lemma
\ref{lemm6.3} with its proof, $v(K[x _{n}])/v(K)$ and $v(K[y
_{n}])/v(K)$ are cyclic groups generated by the cosets
$(1/n)v(\alpha ) + v(K)$ and $(1/n)v(\beta ) + v(K)$,
respectively. Let now $\Delta _{n} \in d(K)$. As shown by Draxl
(cf. \cite[page~216]{Dr1}), then the sum, say $\Sigma (x _{n}, y
_{n}; K)$, of the subgroups $v(K[x _{n}])/v(K)$ and $v(K[y
_{n}])/v(K)$ of $v(\Delta _{n})/v(K)$ is direct. Since $e(\Delta
_{n}/K) = n ^{2}$, and by \cite[Lemma~7.70]{TW}, there is a natural
group isomorphism $(1/n)v(K)/v(K) \cong v(K)/nv(K)$, this implies
$\Sigma (x _{n}, y _{n}; K) = v(\Delta _{n})/v(K)$ and $v(\Delta
_{n})/v(K)$ is isomorphic to the subgroup $\langle v(\alpha ) +
nv(K), v(\beta ) + nv(K)\rangle $ of $v(K)/nv(K)$. Thus the
statements about $v(\Delta _{n})/v(K)$ and the left-to-right
implication in Lemma \ref{lemm6.4} (a) are proved.
\par\vskip0.04truecm
In view of (4.2) (b), it remains to be seen that exp$(\Delta _{n})
= n$ if the subgroup $\langle \alpha K ^{\ast n}, \beta K ^{\ast
n}\rangle $ of $K ^{\ast }/K ^{\ast n}$ has order $n ^{2}$. As the
extension $K[y _{n}]/K$ is totally ramified and cyclic of degree
$n$, and the relation $y _{n}x _{n} = u _{n}x _{n}y _{n}$ gives
rise to a $K$-automorphism $\varphi _{n}$ of $K[y _{n}]$ of order
$n$, such that
\par\vskip0.032truecm\noindent
$\varphi _{n}(y _{n}) = x _{n} ^{-1}y _{n}x _{n} = u _{n}y _{n}$, one
may identify $\Delta _{n}$ with the cyclic $K$-algebra
\par\vskip0.032truecm\noindent
$(K[y \sb n]/K, \varphi _{n}, \alpha )$. Moreover, it follows from
(4.2) (b) and Lemma \ref{lemm6.3} that $\alpha ^{j} \in N(K[y
_{n}]/K)$, where $j \in \mathbb{N}$, if and only if $n \mid j$.
Therefore,
\par\vskip0.032truecm\noindent
\cite[Sect.~15.1, Proposition~b and Corollary~d]{P} show that
exp$(\Delta _{n}) = n$,
\par\vskip0.032truecm\noindent
ind$(\Delta _{n}) = {\rm deg}(\Delta _{n}) = n$, and $\Delta _{n} \in
d(K)$.
\end{proof}
\par
\smallskip
\begin{rema}
\label{rema6.5} The right-to-left implication in the former part
of Lemma \ref{lemm6.4}~{\rm (a)} remains valid if $(K, v)$ is a
valued field, char$(\widehat K) \nmid n$, and $K$ contains a
primitive $n$-th root of unity. Indeed, take a Henselization $(K
^{\prime }, v')$ of $(K, v)$, put $\Delta _{n} ^{\prime } = \Delta
_{n} \otimes _{K} K ^{\prime }$, and observe that $(K ^{\prime },
v')/(K, v)$ is a valued field extension with $\widehat K ^{\prime
} = \widehat K$ and $v'(K ^{\prime }) = v(K)$ (cf.
\cite[Theorem~15.3.5]{E3}). In view of the implication $\Delta
_{n} ^{\prime } \in d(K ^{\prime }) \to \Delta _{n} \in d(K)$,
this allows us to carry out our proof only in the case of $(K, v)$
Henselian. Then $v(K _{\rm ur}) = v(K)$, by Lemma \ref{lemm4.3}
{\rm (b)}, and $\Delta _{n} \otimes _{K} K _{\rm ur}$ is a symbol
$K _{\rm ur}$-algebra, so it follows from
\par\noindent
Lemma \ref{lemm6.4} {\rm (a)} that $\Delta _{n} \otimes _{K} K
_{\rm ur} \in d(K _{\rm ur})$, whence, $\Delta _{n} \in d(K)$ (and
\par\noindent
$v _{K _{\rm ur}}(\Delta _{n} \otimes _{K} K _{\rm ur}) = v(\Delta
_{n})$, by \cite[Theorem~1]{Mo}). Furthermore, by the
Henselization theorem (cf. \cite[Theorem~2]{Mo},
\cite[Proposition~3]{Er}, \cite[Sect.~4.1 and pp.~189-190]{TW}),
the fact that $\Delta _{n} ^{\prime } \in d(K ^{\prime })$
indicates that $v$ extends
\par\vskip0.04truecm\noindent
uniquely to a valuation $v _{\Delta _{n}}$ of $\Delta _{n}$, and
$(\Delta _{n} ^{\prime }, v _{\Delta _{n}'})/(\Delta _{n}, v
_{\Delta _{n}})$ is a valued division ring extension with
$\widehat {\Delta _{n}'} = \widehat \Delta _{n} = \widehat K$ and
$v'(\Delta _{n} ^{\prime }) = v _{\Delta _{n}}(\Delta _{n}) =
\langle (1/n)v(\alpha ), (1/n)v(\beta )\rangle + v(K)$.
\end{rema}
\par
\smallskip Our proof of Proposition \ref{prop2.4} is facilitated by
Lemmas \ref{lemm6.1}, \ref{lemm6.2} and the following lemma, which
make it possible to achieve this goal by valuation-theoretic
arguments relying on the structure theory of finite-dimensional
central simple algebras and avoiding technical complications.
\par
\medskip
\begin{lemm}
\label{lemm6.6} Assume that $(K, v)$ and $\widehat K$ are as in
Lemma \ref{lemm6.4}, fix a finite extension $L/K$ with {\rm
char}$(\widehat K) \nmid [L\colon K]$, and let $D \in d(K)$ be a
tame algebra over $K$. Then $D \otimes _{K} L \in d(L)$ if and
only if $v(D) \cap v(L) = v(K)$; when this holds, $v _{L}(D
\otimes _{K} L) = v(D) + v(L)$ and the group $v _{L}(D \otimes
_{K} L)/v(K)$ is isomorphic to the direct sum $v(D)/v(K) \oplus
v(L)/v(K)$.
\end{lemm}
\par
\smallskip
\begin{proof}
Take an algebra $D _{L} \in d(L)$ so that $[D _{L}] = [D \otimes
_{K} L]$ in Br$(L)$. Then Wedderburn's theorem, applied to $D
\otimes _{K} L$, implies deg$(D _{L}) \mid {\rm deg}(D)$; as $D/K$
is tame and $\widehat K = \widehat K _{\rm sep}$, it follows that
$e(D/K) = [D\colon K]$ and $e(D _{L}/L) = [D _{L}\colon L]$ (i.e.
$D/K$ and $D _{L}/L$ are totally ramified). In addition, deg$(D
_{L}) = {\rm deg}(D)$ if and only if $D \otimes _{K} L \in d(L)$,
i.e. $D \otimes _{K} L = D _{L}$. Note further that $v _{L}(D
_{L})$ is a subgroup of $v(D) + v(L)$ (cf.
\cite[Corollary~6.6]{JW}), whence, $v _{L}(D _{L})/v(L)$ is a
subgroup of $(v(D) + v(L))/v(L)$. The groups $(v(D) + v(L))/v(L)$
and $v(D)/(v(D) \cap v(L))$ are isomorphic, so there is a
surjective homomorphism $\eta _{L/K}$ of $v(D)/v(K)$ upon $(v(D) +
v(L))/v(L)$ ($\eta _{L/K}$ is an isomorphism if and only if $v(D)
\cap v(L) = v(K)$). Since $D/K$ and $D _{L}/L$ are totally
ramified and deg$(D _{L}) \mid {\rm deg}(D)$, these remarks show that
deg$(D _{L}) = {\rm deg}(D)$ if and only if $v(D) \cap v(L) = v(K)$
and
\par\noindent
$v _{L}(D _{L}) = v(D) + v(L)$. They prove the former assertion of
Lemma \ref{lemm6.6}, and when $D \otimes _{K} L \in d(L)$, the
equality $v _{L}(D _{L}) = v(D) + v(L)$. Suppose finally that
$v(D) \cap v(L) = v(K)$. Then the subgroups $v(D)/v(K)$ and
$v(L)/v(K)$ of (the abelian group) $\overline {v(K)}/v(K)$
intersect trivially, i.e. their inner sum is direct and equal to
$(v(D) + v(L))/v(K)$, which completes our proof.
\end{proof}
\par
\medskip
Note here that a field $E$ is called $p$-quasilocal, for some $p
\in \mathbb{P}$, if one of the following two conditions holds:
Brd$_{p}(E) = 0$ or $E(p) = E$; Brd$_{p}(E) \neq 0$, $E(p) \neq E$
and every degree $p$ extension of $E$ in $E(p)$ is embeddable as
an $E$-subalgebra in each $D \in d(E)$ with deg$(D) = p$. By class
field theory, $1$-local fields with finite residue fields are $p$-quasilocal, for all $p
\in \mathbb{P}$. When $E$ is $p$-quasilocal, we have Brd$_{p}(E) \le
1$ if $E(p) \neq E$, as well as in the case where
$p = {\rm char}(E)$ or $E$ contains a primitive $p$-th root of unity
(see \cite[Sect. 3]{Ch3}). Our next lemma characterizes those
strictly Henselian fields of residual characteristic different from
$p$, which possess the $p$-quasilocal property.
\par
\medskip
\begin{lemm}
\label{lemm6.7} Let $(K, v)$ be a strictly Henselian field and let
$p \in \mathbb{P}_{\widehat K}$. Then $K$ is a $p$-quasilocal
field if and only if the group $v(K)/pv(K)$ has order less than $p
^{3}$; this holds if and only if {\rm Brd}$_{p}(K) = 0$ or
$v(K)/pv(K)$ has order $p ^{2}$.
\end{lemm}
\par
\smallskip
\begin{proof}
Evidently, $v(K)/pv(K)$ can be viewed as a vector space over
$\mathbb{Z}/p\mathbb{Z}$ of dimension $r _{p}$, so the assertion
of the lemma can be restated by saying that $K$ is a
$p$-quasilocal field if and only if $r _{p} \le 2$. It follows
from Lemma \ref{lemm6.4} that if $r _{p} \le 1$, then Br$(K) _{p}
= \{0\}$, i.e. Brd$_{p}(K) = 0$; in particular, $K$ is
$p$-quasilocal. Next we consider the case where $r _{p} = 2$.
Using again Lemma \ref{lemm6.4}, we obtain that every $D _{p} \in
d(K)$ with $[D _{p}] \in {\rm Br}(K)_{p}$ is a symbol $K$-algebra,
and if deg$(D _{p}) = p$, then $v(D _{p}) = \{\gamma \in \overline
{v(K)}\colon p\gamma \in v(K)\}$. Together with Lemmas
\ref{lemm4.1} (a) and \ref{lemm6.6}, this implies the
non-existence of a field $L \in I(K(p)/K)$, such that $[L\colon K]
= p$ and $D _{p} \otimes _{K} L \in d(L)$. It is therefore clear
from \cite[Lemma~3.5]{Ch2} (or \cite[Sect. 13.4, Corollary~b]{P})
that degree $p$ extensions of $K$ are embeddable as
$K$-subalgebras in $D _{p}$; in other words, $K$ is
$p$-quasilocal.
\par
Suppose finally that $r _{p} \ge 3$, fix a primitive $p$-th root
of unity $u _{p} \in K$, and take $\pi _{1}, \pi _{2}, \pi _{3}
\in K$ with cosets $(v(\pi _{i}) + pv(K)) \in v(K)/pv(K)\colon i =
1, 2, 3$, linearly independent over $\mathbb{Z}/p\mathbb{Z}$. Put
$\Delta _{p} = K(\pi _{1}, \pi _{2}; u _{p}) _{p}$ and $L = K(\pi
)$, for some $p$-th root $\pi \in K(p)$ of $\pi _{3}$. It follows
from Lemmas \ref{lemm4.1} (a) and \ref{lemm6.4} (a) that $v(\Delta
_{p}) \cap v(L) = v(K)$, whence, by Lemma \ref{lemm6.6}, $\Delta
_{p} \otimes _{K} L \in d(L)$. This observation shows that $K$ is
not $p$-quasilocal, so Lemma \ref{lemm6.7} is proved.
\end{proof}
\par
\medskip
The $p$-quasilocality of $K$ in case $r _{p}(K) = 2$ of Lemma
\ref{lemm6.7} is also implied by \cite[Lemma~3.8]{Ch3} and the
fact that $\mathcal{G}(K(p)/K) \cong \mathbb{Z} _{p} ^{2}$ and
$\mathbb{Z} _{p} ^{2}$ is a Demushkin group (cf.
\cite[Theorem~A.24~(v), Lemma 7.70]{TW} and \cite[Lemma~7]{W}).
\par
\medskip
\begin{rema}
\label{rema6.8} Let $K _{0}$ be a field, $K _{n} = K _{0}((Z
_{1})) \dots ((Z _{n}))$ an $n$-fold iterated Laurent series field
over $K _{0}$, and $\omega _{n}$ the standard $\mathbb{Z}
^{n}$-valued valuation of $K _{n}$ trivial on $K _{0}$, for some
$n \in \mathbb{N}$. Then $(K _{n}, \omega _{n})$ is Henselian
(more precisely, maximally complete) with $\widehat K _{n} = K
_{0}$, so there exist group isomorphisms
$$K _{n} ^{\prime \ast }/K _{n} ^{\prime \ast p} \cong \widehat K
_{n} ^{\prime \ast }/\widehat K _{n} ^{\prime \ast p} \times v
_{n}(K _{n} ^{\prime })/pv _{n}(K _{n} ^{\prime }) \cong K _{0}
^{\ast }/K _{0} ^{\ast p} \times v _{n}(K _{n} ^{\prime })/pv
_{n}(K _{n} ^{\prime }),$$ for any finite extension $K _{n}
^{\prime }/K _{n}$ and each $p \in \mathbb{P}_{K_{0}}$. When $K
_{0} = K _{0,{\rm sep}}$, this
\par\smallskip\noindent
yields $K _{n} ^{\prime \ast }/K _{n} ^{\prime \ast p} \cong v
_{n}(K _{n} ^{\prime })/pv _{n}(K _{n} ^{\prime }) \cong
(\mathbb{Z}/p\mathbb{Z}) ^{n}$.
\end{rema}
\par
\medskip
\section{\bf An analog to Proposition 2.4 over iterated Laurent
series fields}
\par
\medskip
Let $E/E _{0}$ be a field extension, $D$ a central division
$E$-algebra possessing an $E _{0}$-subalgebra $\Delta _{0} \in d(E
_{0})$, and $\Delta $ the $E$-subalgebra of $D$ generated by the
set $\Delta _{0} \cup E$. Then each $\xi \in \Delta $ is
presentable as an $E$-linear combination of elements of any fixed
$E _{0}$-basis of $\Delta _{0}$. This implies $[\Delta \colon E]
\le [\Delta _{0}\colon E _{0}]$ and there exists a surjective
$E$-homomorphism $\eta \colon \Delta _{0} \otimes _{E _{0}} E \to
\Delta $. Note that $\eta $ is an isomorphism: its injectivity,
used in \cite{Cohn}, follows from the assumption that $\Delta _{0}
\in d(E _{0})$ which guarantees that $\Delta _{0} \otimes _{E
_{0}} E \in s(E)$ (cf. \cite[Sect.~12.4, Proposition~b~(ii)]{P}).
These facts are repeatedly used in this section to prove the
following result and thereby to take a major step towards the
proof of Theorem \ref{theo2.1}.
\par
\medskip
\begin{prop}
\label{prop7.1} Let $K _{0}$ be a field with {\rm char}$(K _{0}) =
q$ and $K _{0,{\rm sep}} = K _{0}$, and for a given $p \in
\mathbb{P}_{K _{0}}$, take a sequence $\bar \varepsilon =
\varepsilon _{m}\colon m \in \mathbb{N} \cup \{0\}$, of elements
of $K _{0}$, such that $\varepsilon _{0} = 1$, $\varepsilon _{1}
\neq 1$, and $\varepsilon _{m} ^{p} = \varepsilon _{m-1}$, for
every $m > 0$. For any $n \in \mathbb{N}$, denote by $K _{n}$ the
$2n$-fold iterated Laurent series field $K _{0}((X _{1}))((Y
_{1})) \dots ((X _{n}))((Y _{n}))$, put $L _{n} = K
_{n}(\sqrt[p]{Y _{i}}\colon i = 1, \dots , n) = K _{0}((X
_{1}))((\sqrt[p]{Y _{1}})) \dots ((X _{n}))((\sqrt[p]{Y _{n}}))$,
and let
\par\vskip0.056truecm\noindent
$D _{n} = D _{n}(X _{i}, \sqrt[p]{Y _{i}}\colon i = 1, \dots , n;
\bar \varepsilon , \bar \mu )$ be the $L _{n}$-algebra $\otimes
_{i=1} ^{n} L _{n}(X _{i}, \sqrt[p]{Y _{i}}; \varepsilon _{\mu
_{i}}) _{p ^{\mu _{i}}}$, where $\otimes = \otimes _{L _{n}}$ and
$\bar \mu = (\mu _{1}, \dots , \mu _{n})$ is a fixed $n$-tuple of
integers $\ge 0$. Then $D _{n} \in d(L _{n})$ and $K _{n}$ is the
unique central $K _{n}$-subalgebra of $D _{n}$.
\end{prop}
\par
\smallskip
\begin{proof}
We first show that $D _{n} \in d(L _{n})$. Let $v _{n}$ be the
standard (Henselian) $\mathbb{Z} ^{2n}$-valued valuation of $K
_{n}$ trivial on $K _{0}$, and let $v _{n}'$ be the prolongation
of $v _{n}$ on $L _{n}$. Then $\widehat K _{n} = K _{0} = K
_{0,{\rm sep}}$ and $v _{n}(K _{n}) = \langle v _{n}(X _{i}), v
_{n}(Y _{i})\colon 1 \le i \le n\rangle $, which implies the
system $v _{n}(X _{i}) + pv _{n}(K _{n}), v _{n}(Y _{n}) + pv
_{n}(K _{n})\colon 1 \le i \le n$, is a basis of $v _{n}(K
_{n})/pv _{n}(K _{n})$ (as a vector space over
$\mathbb{Z}/p\mathbb{Z}$). Similarly, the cosets $v _{n}(X _{i}) +
pv _{n}(L _{n}), (1/p)v _{n}(Y _{i}) + pv _{n}(L _{n})\colon 1 \le
i \le n$, form a basis of $v _{n}(L _{n})/pv _{n}(L _{n})$, so it
follows from Lemma \ref{lemm6.4} that $D _{n} \in d(L _{n})$.
\par
We prove by induction on $n$ that $K _{n}$ is the only central $K
_{n}$-subalgebra of $D _{n}$. Suppose first that $n = 1$ (so $D
_{1} = L _{1}(X _{1}, \sqrt[p]{Y _{1}}; \varepsilon _{\mu _{1}}; p
^{\mu _{1}})$). As noted above, $v _{1}$ is strictly Henselian,
and by Remark \ref{rema6.8}, there are isomorphisms $v _{1}(K
_{1})/pv _{1}(K _{1}) \cong (\mathbb{Z}/p\mathbb{Z}) ^{2} \cong K
_{1} ^{\ast }/K _{1} ^{\ast p}$. It is therefore clear from Lemmas
\ref{lemm6.4} and \ref{lemm6.7} that $K _{1}$ is a $p$-quasilocal
field (with Br$(K _{1}) _{p} \cong \mathbb{Q}_{p}/\mathbb{Z}_{p}
\neq \{0\}$, see \cite[Theorem~3.1]{Ch3}). Moreover, by
\cite[Theorem~4.1]{Ch3}, ind$(B \otimes _{K _{1}} L _{1}) = {\rm
deg}(B)/p$, for every $B \in d(K _{1})$ with $[B] \in {\rm Br}(K
_{1}) _{p}$ and $[B] \neq 0$. On the other hand, the remarks at
the beginning of this section show that if $D _{1}$ has a central
$K _{1}$-subalgebra $\Sigma _{1} \neq K _{1}$, then $\Sigma _{1}
\in d(K _{1})$ and $\Sigma _{1} \otimes _{K _{1}} L _{1}$ is
isomorphic to the $L _{1}$-subalgebra $\Sigma _{1} ^{\prime } =
\langle \Sigma _{1} \cup L _{1}\rangle $ of $D _{1}$. This
requires that $\Sigma _{1} ^{\prime } \in d(L _{1})$ and
ind$(\Sigma _{1} ^{\prime }) = {\rm ind}(\Sigma _{1})$, which is
impossible, since $K _{1}$ is $p$-quasilocal, $L _{1} \neq K _{1}$
and $L _{1} \in I(K _{1}(p)/K _{1})$. Thus Proposition
\ref{prop7.1} is proved in the case where $n = 1$.
\par
\medskip
For the rest of the proof, we need the following lemma.
\par
\medskip
\begin{lemm}
\label{lemm7.2} Assume that $(K, v)$ is a Henselian field and $D =
\Delta \otimes _{K} T$, where $\Delta $ and $T \in d(K)$, $\Delta
/K$ is inertial, $T/K$ is totally ramified and {\rm deg}$(T)$ is
not divisible by {\rm char}$(\widehat K)$. Then $D \in d(K)$,
$\widehat D = \widehat \Delta $, and $v(D) = v(T)$.
\end{lemm}
\par
\smallskip
\begin{proof}
This follows from \cite[Theorem~1]{Mo} as well as from
\cite[Corollary~6.8]{JW}.
\end{proof}
\par
\medskip
Henceforth, we suppose that $n \ge 2$ and the latter part of the
statement of Proposition \ref{prop7.1} holds if $n$ is replaced by
any integer $\nu $ with $1 \le \nu < n$. We prove this part for
$n$ by assuming the opposite, and regarding the $L _{n}$-algebra
$D _{n} = D _{n}(X _{i}, \sqrt[p]{Y _{i}}\colon i = 1, \dots , n;
\bar \varepsilon ; \bar \mu )$ as a counter-example to the
considered assertion, for which deg$(D _{n}) = \prod _{i=1} ^{n} p
^{\mu _{i}}$ is minimal. Taking $v _{n}$ and $v _{n}'$ as above,
let $w _{n}$ be the $\mathbb{Z} ^{2}$-valued valuation of $K _{n}$
trivial on the $(2n - 2)$-fold
\par\smallskip\noindent
Laurent series field $K _{n-1} = K _{0}((X _{1}))((Y _{1})) \dots
((X _{n-1}))((Y _{n-1}))$, and let $w _{n}'$ be the valuation of
$L _{n}$ extending $w _{n}$. Put
$$\widetilde D _{n} = D _{n}(X _{i}, \sqrt[p]{Y _{i}}\colon i = 1,
\dots , n; \bar \varepsilon ; \tilde \mu ) := \otimes _{i=1} ^{n}
L _{n}(X _{i}, \sqrt[p]{Y _{i}}; \varepsilon _{\tilde \mu _{i}})
_{p ^{\tilde \mu _{i}}},$$ where $\otimes = \otimes _{L _{n}}$ and
$\tilde \mu = (\tilde \mu _{1}, \dots , \tilde \mu _{n})$ is
defined as follows: $\tilde \mu _{i} = \mu _{i} - 1$ if
\par\vskip0.04truecm\noindent
$\mu _{i} > 0$; $\tilde \mu _{i} = 0$ if $\mu _{i} = 0$. It
follows from Kummer theory and the theory of cyclic algebras over
arbitrary fields that $[\widetilde D _{n}] = p[D _{n}]$ (cf.
\cite[Sect. 15.1, Corollary~b]{P}). Considering $(K _{n}, v _{n})$
and applying Lemma \ref{lemm6.2} and \cite[Proposition~6.9]{JW},
one obtains further that $\widetilde D _{n} \in d(L _{n})$. Recall
here that we have already seen that if $\Sigma _{n} \in d(K _{n})$
is a $K _{n}$-subalgebra of $D _{n}$, and $\Sigma _{n} ^{\prime }$
is the $L _{n}$-subalgebra $\langle \Sigma _{n} \cup L _{n}\rangle
$ of $D _{n}$, then $\Sigma _{n} ^{\prime } \cong \Sigma _{n}
\otimes _{K _{n}} L _{n}$ over $L _{n}$. Identifying $\Sigma _{n}
\otimes _{K _{n}} L _{n}$ with $\Sigma _{n} ^{\prime }$, one
deduces from the double centralizer theorem that
\par\vskip0.04truecm\noindent
$D _{n} \cong (\Sigma _{n} \otimes _{K _{n}} L _{n}) \otimes _{L
_{n}} \Psi _{n}$, where $\Psi _{n} = C _{D
_{n}}(\Sigma _{n} \otimes _{K _{n}} L _{n})$. Now choose
\par\vskip0.04truecm\noindent
algebras $\widetilde \Sigma _{n} \in d(K _{n})$ and $\widetilde \Psi
_{n} \in d(L _{n})$ so that $[\widetilde \Sigma _{n}] = p[\Sigma
_{n}]$ and $[\widetilde \Psi _{n}] = p[\Psi _{n}]$, in Br$(K _{n})$
and Br$(L _{n})$, respectively. Taking into account that
\par\vskip0.04truecm\noindent
$\Sigma _{n} \otimes _{K _{n}} L _{n} \in d(L _{n})$, and $v
_{n}(\widetilde \Sigma _{n}) = pv _{n}(\Sigma _{n}) + v _{n}(K
_{n}) \subset v _{n}(\Sigma _{n})$ (by
\par\vskip0.04truecm\noindent
\cite[Proposition~6.9]{JW}), and applying Lemma \ref{lemm6.6}, one
obtains consecutively
\par\smallskip\noindent
$${\rm that} \ v _{n}(\widetilde \Sigma _{n}) \cap v _{n}(L _{n}) = v
_{n}(\Sigma _{n}) \cap v _{n}(L _{n}) = v _{n}(K _{n}) \ {\rm and}
\ \widetilde \Sigma _{n} \otimes _{K _{n}} L _{n} \in d(L _{n}).$$
\par\smallskip\noindent
Similarly, it follows from \cite[Lemma~3]{Dr1} and the $L
_{n}$-isomorphism
\par\smallskip\noindent
$D _{n} \cong (\Sigma _{n} \otimes _{K _{n}} L _{n}) \otimes _{L
_{n}} \Psi _{n}$ that $v _{n}'(\Sigma _{n} \otimes _{K _{n}} L
_{n}) \cap v _{n}'(\Psi _{n}) = v _{n}'(L _{n}) (:= v _{n}(L
_{n}))$.
\par\smallskip\noindent
In view of \cite[Proposition~6.9]{JW}, this implies
\par\vskip0.04truecm\noindent
$v _{n}'(\widetilde \Sigma _{n} \otimes _{K _{n}} L _{n}) \cap v
_{n}'(\widetilde \Psi _{n}) = v _{n}'(L _{n})$, so it can be
deduced from \cite[Lemma~3]{Dr1} that $\widetilde D _{n} \cong
(\widetilde \Sigma _{n} \otimes _{K _{n}} L _{n}) \otimes _{L
_{n}} \widetilde \Psi _{n}$ as $L _{n}$-algebras. Thus it turns
out that $\widetilde \Sigma _{n}$ embeds in $\widetilde D _{n}$ as
a central $K _{n}$-subalgebra, and by the choice of $D _{n}$, this
means that $\widetilde \Sigma _{n} = K _{n}$ and exp$(\Sigma _{n})
= p$. Now it follows from \cite[Theorem~1]{Dr1} that $\Sigma _{n}$
decomposes into a tensor product of symbol $K _{n}$-algebras of
degree $p$, which allows us to assume for the rest of the proof of
Proposition \ref{prop7.1} that $\Sigma _{n}$ is chosen to be such
a symbol $K _{n}$-algebra.
\par\vskip0.05truecm
By (4.2) (b), we have $K _{n} ^{\ast }/K _{n} ^{\ast p} \cong v
_{n}(K)/pv _{n}(K)$, so one may view $K _{n} ^{\ast }/K _{n}
^{\ast p}$
\par\vskip0.04truecm\noindent
as a $2n$-dimensional vector space over $(\mathbb{Z}/p\mathbb{Z})$
with a basis formed by cosets
\par\vskip0.04truecm\noindent
$X _{\nu }.K _{n} ^{\ast p}$, $Y _{\nu }.K _{n} ^{\ast p}$, $\nu =
1, \dots , n$. Observe also that every symbol $K _{n}$-subalgebra
of $\Sigma _{n} \otimes _{K _{n}} L _{n}$ of degree $p$ is
embeddable in $D _{n}$, $L _{n} ^{\ast p}$ contains the elements
\par\vskip0.04truecm\noindent
$Y _{1}, \dots , Y _{n}$, and each $\lambda \in K _{n} ^{\ast }
\setminus K _{n} ^{\ast p}$ belongs to the norm group $N(K
_{n}(\sqrt[p]{\lambda })/K _{n})$.
\par\vskip0.04truecm\noindent
Therefore, $\Sigma _{n}$ can be chosen so as to be isomorphic to
$K _{n}(a(\widetilde X), b(\widetilde X).X ^{j}; \varepsilon _{1})
_{p}$ over $K _{n}$, for some monic monomials $a(\widetilde X)$
and $b(\widetilde X) \in K _{0}[\widetilde X]$, where $j \in \{0,
1\}$ and $\widetilde X$ is the $n - 1$-tuple $(X _{i'}\colon 1 \le
i' \le n - 1)$.
\par\vskip0.04truecm
In order to prove that the existence of $\Sigma _{n}$ is
impossible, it is convenient, in the rest of our proof, to
consider finite extensions of $K _{n}$ and finite-dimensional
division $K _{n}$-algebras (including $\Sigma _{n}$ and other
subalgebras of $D _{n}$) with their valuations extending $w _{n}$;
to simplify notation, we write $w _{n}'$ for the prolongation of
$w _{n}$ on $L _{n}$. The algebra $\Sigma _{n}$ is clearly
inertial or NSR over $K _{n}$ depending on whether $j = 0$ or $j =
1$; note also that $a(\tilde X) \in K _{n-1} ^{\ast } \setminus L
_{n} ^{\ast p}$,
\par\vskip0.04truecm\noindent
and $\Sigma _{n} \otimes _{K _{n}} L _{n} \cong \Sigma _{n}
^{\prime } \cong L _{n}(a(\widetilde X), b(\widetilde X);
\varepsilon _{1}) _{p}$ as $L _{n}$-algebras. Put
\par\smallskip\noindent
$L _{n-1} = K _{n-1}(\sqrt[p]{Y _{i}}\colon i = 1, \dots , n -
1)$, $L _{n-1} ^{\prime } = K _{n}(\sqrt[p]{Y _{i}}\colon i = 1,
\dots , n - 1)$,
\par\smallskip\noindent
$D _{n-1} = \otimes _{i=1} ^{n-1} L _{n-1}(X _{i}, \sqrt[p]{Y
_{i}}; \varepsilon _{\mu _{i}}) _{p ^{\mu _{i}}}$, where $\otimes
= \otimes _{L _{n-1}}$, and
\par\smallskip\noindent
$D _{n-1} ^{\prime } = \otimes _{i=1} ^{n-1} L _{n}(X _{i},
\sqrt[p]{Y _{i}}; \varepsilon _{\mu _{i}}) _{p ^{\mu _{i}}}$,
where $\otimes = \otimes _{L _{n}}$.
\par\smallskip\noindent
It is easily verified that $L _{n} = L _{n-1}((X
_{n}))((\sqrt[p]{Y _{n}}))$, $\widehat L _{n} = L _{n-1}$,
\par\vskip0.04truecm\noindent
$D _{n-1} \in d(L _{n-1})$, and the $L _{n}$-algebras $D _{n-1}
^{\prime } = D _{n-1}((X _{n}))((\sqrt[p]{Y _{n}}))$ and
\par\vskip0.04truecm\noindent
$D _{n-1} \otimes _{L _{n-1}} L _{n}$ are isomorphic. Moreover,
Lemma \ref{lemm4.2}~(a) and
\par\vskip0.04truecm\noindent
\cite[Theorem~2.8]{JW} imply that $D _{n-1} ^{\prime } \in d(L
_{n})$, $D _{n-1} ^{\prime }$ is a maximal inertial $L
_{n}$-subalgebra of $D _{n}$ (its uniqueness, up-to conjugacy in
$D _{n}$, follows from \cite[Theorem~2.8]{JW} and Skolem-Noether's
theorem), and $\widehat D _{n-1} ^{\prime } = \widehat D _{n}
\cong D _{n-1}$ as an $L _{n-1}$-algebra. In this setting, we
complete the proof of Proposition \ref{prop7.1} in the following
three steps.
\par\smallskip\noindent
Step 1. We show that $D _{n}$ does not possess an inertial symbol
$K _{n}$-subalgebra of degree $p$. Assuming the opposite, we may
suppose that $j = 0$, i.e. $\Sigma _{n}/K _{n}$ is inertial. As $w
_{n}(a(\tilde X)) = w _{n}(b(\tilde X)) = 0$, this means that
there are
\par\smallskip\noindent
isomorphisms $\Sigma _{n} \cong K _{n}(a(\tilde X), b(\tilde X);
\varepsilon _{1}) _{p}$ and $\Sigma _{n} ^{\prime } \cong L
_{n}(a(\tilde X), b(\tilde X); \varepsilon _{1}) _{p}$ over
\par\smallskip\noindent
$K _{n}$ and $L _{n}$, respectively; in particular, $\Sigma _{n}
^{\prime }/L _{n}$ is inertial. Therefore, by Lemma \ref{lemm6.1},
$\Sigma _{n} ^{\prime }$ embeds in $D _{n-1} ^{\prime }$ as an $L
_{n}$-subalgebra, and by \cite[Theorem~2.8]{JW}, so does $\widehat
\Sigma _{n} ^{\prime }$ in $\widehat D _{n-1} ^{\prime } = D
_{n-1}$ as an $L _{n-1}$-subalgebra. Since $\widehat \Sigma _{n}
^{\prime }$,
\par\vskip0.052truecm\noindent
$L _{n-1}(a(\tilde X)), b(\tilde X); \varepsilon _{1}) _{p}$ and
$K _{n-1}(a(\tilde X), b(\tilde X); \varepsilon _{1}) _{p} \otimes
_{K _{n-1}} L _{n-1}$ are
\par\vskip0.052truecm\noindent
isomorphic $L _{n-1}$-algebras (see \cite[Theorem~3.1]{JW}), this
contradicts the inductive hypothesis that the assertion of
Proposition \ref{prop7.1} holds if $n$ is replaced by $n - 1$. The
encountered contradiction proves the non-existence of degree $p$
inertial symbol $K _{n}$-subalgebras of $D _{n}$; in particular,
it follows that $\Sigma _{n}/K _{n}$ is not inertial, which yields
$j = 1$. It is also clear that $\mu _{n} > 0$, i.e. $L _{n}(X
_{i}, \sqrt[p]{Y _{i}}; \varepsilon _{\mu _{i}}) _{p ^{\mu _{i}}}
\neq L _{n}$.
\par\smallskip\noindent
Step 2. We prove that the $L _{n}$-algebra $B _{n} = L _{n}(\beta
_{n-1}, X _{n}; \varepsilon _{1}) _{p}$ does not embed in $D _{n}$ as
an $L _{n}$-subalgebra, for any $\beta _{n-1} \in K
_{n-1} ^{\ast } \setminus L _{n-1} ^{\ast p}$ (this statement is a
part of (7.1) (b)). Our argument relies on Kummer theory and
\cite[Sect. 15.1, Proposition~b]{P}, which imply that $B _{n}$ and
its opposite $L _{n}$-algebra $B _{n} ^{\rm op}$ are cyclic,
and $B _{n} ^{\rm op}$ is isomorphic to $L _{n}(\beta _{n-1}, X
_{n} ^{-1}; \varepsilon _{1}) _{p}$ and
\par\vskip0.04truecm\noindent
$L _{n}(X _{n}, \beta _{n-1}; \varepsilon _{1}) _{p}$. Identifying
$B _{n} ^{\rm op}$ with $L _{n}(X _{n}, \beta _{n-1}; \varepsilon
_{1}) _{p}$, and putting
\par\vskip0.04truecm\noindent
$\Lambda _{n} = L _{n}(X _{n}, \beta _{n-1} ^{p ^{\tilde \mu
_{n}}}\sqrt[p]{Y _{n}}; \varepsilon _{\mu _{n}}) _{p ^{\mu _{n}}}$,
we show that the underlying division
\par\smallskip\noindent
$L _{n}$-algebras $\mathcal{D} _{n}$ and $\mathcal{T} _{n}$ of $D
_{n} \otimes _{L _{n}} B _{n} ^{\rm op}$ and $T _{n} \otimes _{L
_{n}} B _{n} ^{\rm op}$, respectively, satisfy the following:
\par
\medskip\noindent
(7.1) (a) $\mathcal{T} _{n}$ and $\Lambda _{n}$ are isomorphic; in
addition, $\mathcal{T} _{n}/L _{n}$ is totally ramified with $w
_{n}'(\mathcal{T} _{n}) = w _{n}'(T _{n}) = w _{n}'(D _{n})$;
\par\vskip0.04truecm
(b) The $L _{n}$-algebras $T _{n} \otimes _{L _{n}} B _{n} ^{\rm
op}$ and $D _{n} \otimes _{L _{n}} B _{n} ^{\rm op}$ are
isomorphic to the matrix rings $M _{p}(\Lambda _{n})$ and $M
_{p}(\mathcal{D} _{n})$, respectively; therefore, $B _{n}$ is not
isomorphic to any $L _{n}$-subalgebra of $D _{n}$;
\par\vskip0.04truecm
(c) $\mathcal{D} _{n}$ is isomorphic to the $L _{n}$-algebra $D
_{n-1} ^{\prime } \otimes _{L _{n}} \Lambda _{n}$; in particular,
$D _{n-1} ^{\prime }$ is isomorphic to a maximal inertial
subalgebra of $\mathcal{D} _{n}$.
\par
\medskip\noindent
It follows from Kummer theory and \cite[Sect. 15.1,
Corollary~b]{P} that
\par\noindent
$[L _{n}(X _{n}, \beta _{n-1}; \varepsilon _{1}) _{p}] = [L _{n}(X
_{n}, \beta _{n-1} ^{p ^{\tilde \mu _{n}}}; \varepsilon _{\mu
_{n}}) _{p ^{\mu _{n}}}]$. This enables one to obtain
\par\vskip0.054truecm\noindent
from the description of the relative Brauer groups of cyclic field
extensions (cf. \cite[Sect. 15.1, Proposition~b]{P}) that
$[\mathcal{T} _{n}] = [\Lambda _{n}]$. Note further that $\Lambda
_{n} \in d(L _{n})$ and $\Lambda _{n}/L _{n}$ is totally ramified. As
explained in Remark \ref{rema6.5}, to prove this fact it suffices to
observe that $w _{n}'(\sqrt[p]{Y _{n}}) = (1/p)w _{n}(Y _{n})$, $w
_{n}'(\beta _{n-1}) = 0$,
\par\smallskip\noindent
$w _{n}(L _{n,{\rm ur}}) = w _{n}(L _{n}) = \langle w _{n}(X
_{n}), (1/p)w _{n}(Y _{n})\rangle \cong \mathbb{Z} ^{2}$ and the
cosets
\par\smallskip\noindent
$w _{n}(X _{n}) + p ^{\mu _{n}}w _{n}(L _{n})$, $w _{n}'(\beta
_{n-1} ^{p ^{\tilde \mu }}.\sqrt[p]{Y _{n}}) + p ^{\mu _{n}}w
_{n}(L _{n})$ generate the group
\par\smallskip\noindent
$w _{n}(L _{n})/p ^{\mu _{n}}w _{n}(L _{n}) \cong (\mathbb{Z}/p
^{\mu _{n}}\mathbb{Z}) ^{2}$. These facts show that the symbol
algebra
\par\vskip0.063truecm\noindent
$L _{n,{\rm ur}}(X _{n}, \beta _{n-1} ^{p ^{\tilde \mu
_{n}}}\sqrt[p]{Y _{n}}; \varepsilon _{\mu _{n}}) _{p ^{\mu _{n}}}
\cong \Lambda _{n} \otimes _{L _{n}} L _{n,{\rm ur}}$ lies in $d(L
_{n,{\rm ur}})$, which proves
\par\vskip0.063truecm\noindent
the claimed properties of $\Lambda _{n}$. As $\mathcal{T} _{n} \in d(L
_{n})$ and $[\Lambda _{n}] = [\mathcal{T}_{n}]$, the assertion
that $\Lambda _{n} \cong \mathcal{T}_{n}$ follows from these
properties and Wedderburn's theorem; also, it becomes clear that
$w _{n}'(\Lambda _{n}) = w _{n}'(\mathcal{T} _{n})$. The
equalities $e(T _{n}/L _{n}) = [T _{n}\colon L _{n}]$
\par\vskip0.063truecm\noindent
and $w _{n}'(T _{n}) = w _{n}'(\mathcal{T} _{n})$ are proved in the
same way as the corresponding assertions about $\Lambda _{n}$,
whereas the fact that $w _{n}'(\mathcal{D} _{n}) = w
_{n}'(\mathcal{T} _{n})$ and
\par\vskip0.048truecm\noindent
$\mathcal{D} _{n} \cong D _{n-1} ^{\prime } \otimes _{L _{n}} \Lambda
_{n}$ as $L _{n}$-algebras can be viewed as a special case of
\par\vskip0.048truecm\noindent
Lemma \ref{lemm7.2}. Thus statements (7.1) (a) and (c) are proved.
Since
\par\vskip0.048truecm\noindent
$M _{\varkappa }(A) \cong A \otimes _{L _{n}} M _{\varkappa }(L
_{n})$ over $L _{n}$ (which yields deg$(M _{\varkappa }(A)) =
\varkappa .{\rm deg}(A)$, for any $A \in d(L _{n})$, $\varkappa \in
\mathbb{N}$), these observations and the equalities
\par\smallskip\noindent
$[T _{n} \otimes B _{n} ^{\rm op}] = [\mathcal{T} _{n}]$, deg$(T _{n}
\otimes _{L _{n}} B _{n} ^{\rm op}) = p.{\rm deg}(\mathcal{T} _{n})$,
and $[D _{n} \otimes _{L _{n}} B _{n} ^{\rm op}] = [\mathcal{D}
_{n}]$, deg$(D _{n} \otimes _{L _{n}} B _{n} ^{\rm op}) = p.{\rm
deg}(\mathcal{D} _{n})$, enable one to deduce
the former part of (7.1) (b) from Wedderburn's theorem. The latter
part of (7.1) (b) follows from the former one and
\cite[Lemma~3.5]{Ch2}.
\par\smallskip
Step 3. Retaining the assumption that $\Sigma _{n} ^{\prime } = L
_{n}(a(\tilde X), b(\tilde X)X _{n}; \varepsilon _{1}) _{p}$ is an
\par\vskip0.056truecm\noindent
$L _{n}$-subalgebra of $D _{n}$, we show that the $L _{n}$-algebra
$\Delta _{n} ^{\prime } = L _{n}(a(\tilde X), b(\tilde X);
\varepsilon _{1}) _{p}$
\par\vskip0.056truecm\noindent
is isomorphic to a subalgebra of $D _{n-1} ^{\prime }$. Our proof
is based on the fact that
\par\vskip0.04truecm\noindent
$B _{n} ^{\prime {\rm op}} \cong L _{n}(a(\tilde X), X _{n} ^{-1};
\varepsilon _{1}) _{p}$ over $L _{n}$, where $B _{n} ^{\prime } =
L _{n}(a(\tilde X), X _{n}; \varepsilon _{1}) _{p}$; this
\par\smallskip\noindent
implies $[\Delta _{n} ^{\prime }] = [\Sigma _{n} ^{\prime }] - [B
_{n} ^{\prime }]$ (in Br$(L _{n})$). Note also that (by the
triviality of
\par\smallskip\noindent
the valuation of $K _{n-1}$ induced by $w _{n-1}$) $w
_{n}'(a(\tilde X)) = w _{n}(a(\tilde X)) = 0$, so  the
\par\smallskip\noindent
latter part of (7.1) (b) shows that $B _{n} ^{\prime }$ does not
embed in $D _{n}$ as an
\par\noindent
$L _{n}$-subalgebra. Since $\Sigma _{n} ^{\prime }$ is an $L
_{n}$-subalgebra of $D _{n}$, this proves that $\Sigma _{n}
^{\prime }$ and $B _{n} ^{\prime }$ are not $L _{n}$-isomorphic.
Therefore, $[\Delta _{n} ^{\prime }] \neq 0$, and it follows from
Wedderburn's theorem, the primality of $p$, and the equality
deg$(\Delta _{n} ^{\prime }) = p$
\par\vskip0.054truecm\noindent
that $\Delta _{n} ^{\prime } \in d(L _{n})$ and $\Delta _{n} \in
d(K _{n})$. As $w _{n}(a(\tilde X)) = w _{n}(b(\tilde X)) = 0$, it
is now easy to see that $\Delta _{n}/K _{n}$ and $\Delta _{n}
^{\prime }/L _{n}$ are inertial. In addition, it becomes
\par\vskip0.054truecm\noindent
clear that $\Sigma _{n} ^{\prime } \otimes _{L _{n}} B _{n}
^{\prime {\rm op}} \cong M _{p}(\Delta _{n} ^{\prime })$ as $L
_{n}$-algebras. At the same time, by the
\par\vskip0.054truecm\noindent
double centralizer theorem (as stated in \cite[Sect.~12.7]{P}),
\par\vskip0.054truecm\noindent
$D _{n} = \Sigma _{n} ^{\prime } \otimes _{L _{n}} C _{D
_{n}}(\Sigma _{n} ^{\prime })$. Applying (7.1) to $B _{n} ^{\prime
}$, one proves that
\par\vskip0.18truecm
$D _{n} \otimes _{L _{n}} B _{n} ^{\prime {\rm op}} \cong M
_{p}(\mathcal{D} _{n}) \cong M _{p}(\Delta _{n} ^{\prime } \otimes
_{L _{n}} C _{D _{n}}(\Sigma _{n} ^{\prime }))$ as $L
_{n}$-algebras.
\par\vskip0.18truecm\noindent
These isomorphisms and Wedderburn's theorem, together with the
fact that
\par\vskip0.04truecm\noindent
$\Delta _{n} ^{\prime } \otimes _{L _{n}} C _{D _{n}}(\Sigma _{n}
^{\prime }) \in s(L _{n})$, and with well-known isomorphisms of
tensor
\par\vskip0.04truecm\noindent
products of full matrix $L _{n}$-algebras (see \cite[Sect.~9.3,
Corollary~b]{P}), show that
\par\vskip0.14truecm
$\Delta _{n} ^{\prime } \otimes _{L _{n}} C _{D _{n}}(\Sigma
_{n} ^{\prime }) \in d(L _{n})$ and, more precisely, $\Delta _{n}
^{\prime } \otimes _{L _{n}} C _{D _{n}}(\Sigma _{n} ^{\prime })
\cong \mathcal{D} _{n}$.
\par\vskip0.14truecm\noindent
Observing also that, by (7.1) (c), $D _{n-1} ^{\prime }$ is
isomorphic to any maximal inertial central $L _{n}$-subalgebra of
$\mathcal{D} _{n}$, one gets from Lemma \ref{lemm6.1} that $\Delta
_{n} ^{\prime }$ is isomorphic to an $L _{n}$-subalgebra of $D
_{n-1} ^{\prime }$. As $\widehat L _{n} = L _{n-1}$ and $\widehat
D _{n-1} ^{\prime } = D _{n-1}$, it is now
\par\noindent
easy to see that the residue division algebra $\Delta _{n-1}
^{\prime } = L _{n-1}(a(\tilde X), b(\tilde X); \varepsilon _{1})
_{p}$
\par\vskip0.04truecm\noindent
of $(\Delta _{n} ^{\prime }, w _{n,\Delta _{n}'}')$ is embeddable
in $D _{n-1}$ an $L _{n-1}$-subalgebra. Since
\par\vskip0.04truecm\noindent
$\Delta _{n-1} ^{\prime } \cong \Delta _{n-1} \otimes _{K _{n-1}}
L _{n-1}$ as $L _{n-1}$-algebras, where $\Delta _{n-1}$ is the
symbol
\par\vskip0.048truecm\noindent
algebra $K _{n-1}(a(\tilde X), b(\tilde X); \varepsilon _{1}) _{p}
\in d(K _{n-1})$, this leads to the conclusion that
\par\vskip0.04truecm\noindent
$\Delta _{n-1}$ lies in $d(K _{n-1})$ and embeds in $D _{n-1}$ as
a $K _{n-1}$-subalgebra. Our conclusion, however, contradicts the
inductive hypothesis, so it follows that $D _{n}$ does not possess
a central $K _{n}$-subalgebra $\Sigma _{n} \neq K _{n}$, as
required by Proposition \ref{prop7.1}.
\end{proof}
\par
\smallskip
\section{\bf Proof of Theorem \ref{theo2.1}}
\par
\smallskip
It has been clarified in Section 2 that Proposition \ref{prop7.1}
implies Proposition \ref{prop2.4}, so we are prepared to prove
Theorem \ref{theo2.1} in the special case where $K _{0}$ is an
algebraically closed field, $K = K _{\infty }$ and trd$(K/K _{0})$
is countable. We begin with a lemma which ensures, for any fixed
$p \in \mathbb{P}_{K_{0}}$, the existence of division $K
_{n}$-algebras $R _{n}$, $n \in \mathbb{N}$, satisfying the
conditions stated in Section 2.
\par
\medskip
\begin{lemm}
\label{lemm8.1} Let $K$ be a purely transcendental extension of a
field $K _{0} = K _{0,{\rm sep}}$ with {\rm trd}$(K/K _{0})$
countable, and let $p \in \mathbb{P}_{K _{0}}$ and $\varepsilon
_{n} = \varepsilon _{n}(p)$, $n \in \mathbb{N}$, be elements of $K
_{0}$, such that $\varepsilon _{1} \neq 1 = \varepsilon _{1}
^{p}$, and $\varepsilon _{n+1} ^{p} = \varepsilon _{n}$, for every
$n$. Suppose that $\Sigma _{\infty } = \{X _{n}, Y _{n}\colon n
\in \mathbb{N}\}$ is a transcendence basis and a generating set of
$K/K _{0}$, and for each $n \in \mathbb{N}$, put $\Sigma _{n} =
\{X _{i}, Y _{i}\colon i = 1, \dots , n\}$, $K _{n} = K
_{0}(\Sigma _{n})$, and denote by $\widetilde \Sigma _{n}$ the
subgroup of $K _{n} ^{\ast }$ generated by $\Sigma _{n}$. Then
there exist division $K _{n}$-algebras $R _{n} = R _{n}\{p\}\colon
n \in \mathbb{N}$, and generating sets $\Sigma _{n} ^{\prime } =
\{X _{n,j}, Y _{n,j}\colon j = 1, \dots , n\}$ of $\widetilde
\Sigma _{n}$, $n \in \mathbb{N}$, with the following properties:
\par
{\rm (a)} For each $n \in \mathbb{N}$, $R _{n}$ is a $K
_{n}$-subalgebra of $R _{n+1}$, and the centre $Z _{n} = Z(R
_{n})$ is an extension of $K _{n}$ obtained by adjunction of
$p$-th roots $Y _{n,1} ^{\prime }, \dots , Y _{n,n} ^{\prime }$ of
$Y _{n,1}, \dots , Y _{n,n}$, respectively; in addition, $Z _{n}
\cap Z _{2n} = K _{n}$ and $K _{n}$ is the unique central $K
_{n}$-subalgebra of $R _{n}$;
\par
{\rm (b)} $R _{n}$ is isomorphic as a $Z _{n}$-algebra to $\otimes
_{u=1} ^{n} S _{n,u}$, where $\otimes = \otimes _{Z _{n}}$ and $S
_{n,u} = Z _{n}(X _{n,u}, Y _{n,u} ^{\prime }; \varepsilon _{\mu
(n, u)}) _{p ^{\mu (n, u)}}$, $\mu (n, u)$ being a positive
integer, for each index $u$; in particular, $[R _{n}\colon K
_{n}]$ is a $p$-power.
\end{lemm}
\par
\smallskip
\begin{proof}
The sequence $R _{n}\colon n \in \mathbb{N}$, is defined
inductively. First we put $$X _{1,1} = X _{1}, Y _{1,1} = Y _{1},
Z _{1} = K _{1}(Y _{1,1} ^{\prime }), R _{1} = S _{1,1} = Z _{1}(X
_{1,1}, Y _{1,1} ^{\prime }; \varepsilon _{1}) _{p}$$
$$= Z _{1} \langle \xi _{1,1}, \eta _{1,1}\colon \xi _{1,1} ^{p} = X
_{1,1}, \eta _{1,1} ^{p} = Y _{1,1} ^{\prime }, \eta _{1,1}\xi
_{1,1} = \varepsilon _{1}\xi _{1,1}\eta _{1,1}\rangle .$$ This
yields $\mu (1, 1) = 1$.
\par
Suppose now that the division $K _{0}$-algebra $R _{k}$ has
already been defined, for some $k \in \mathbb{N}$, and it has the
following properties:
\par\medskip\noindent
(8.1) (i) $Z(R _{k}) = Z _{k} = K _{k}(Y _{k,i} ^{\prime }\colon i
= 1, \dots , k)$, where $Y _{k,i} ^{\prime p} = Y _{k,i}\colon 1
\le i \le k$, for some set $\{Y _{k,i}\colon i = 1, \dots , k\} :=
\Sigma _{k} ^{\prime }$ of generators of $\widetilde \Sigma _{k}$.
\par
(ii) $R _{k}$ is isomorphic as a $Z _{k}$-algebra to $\otimes
_{u=1} ^{k} S _{k,u}$, where $\otimes = \otimes _{Z _{k}}$ and
\par\smallskip\noindent
$S _{k,u} = Z _{k} \langle \xi _{k,u}, \eta _{k,u}\colon \eta
_{k,u}\xi _{k,u} = \varepsilon _{\mu (k,u)}\xi _{k,u}\eta _{k,u},
\xi _{k,u} ^{p^{\mu (k,u)}} = X _{k,u}, \eta _{k,u} ^{p^{\mu
(k,u)}} = Y _{k,u} ^{\prime }\rangle $
\par\medskip\noindent
$= Z _{k}(X _{k,u}, Y _{k,u} ^{\prime }; \varepsilon _{\mu (k,
u)}) _{p ^{\mu (k, u)}}$, for some integers $\mu (k, u) > 0$ and
each index $u$.
\par\medskip\noindent
The definition of $R _{\tilde k}$, where $\tilde k = k + 1$,
relies on the fact that there exists an automorphism $\varphi
_{k}$ of $S _{k,1}$ as a $K _{k}$-algebra, such that $\varphi
_{k}(\eta _{k,1}) = \varepsilon _{1+\mu (k, 1)}\eta _{k,1}$,
\par\vskip0.054truecm\noindent
$\varphi _{k}(\xi _{k,1}) = \xi _{k,1}$, and in case $k \ge 2$,
$\varphi _{k}(Y _{k,i} ^{\prime }) = Y _{k,i} ^{\prime }\colon 2
\le i \le k$. Since the field
\par\smallskip\noindent
$Z _{k}$ is isomorphic over $K _{k}$ to $K _{k}(Y _{k,1} ^{\prime
}) \otimes _{K _{k}} K _{k}(Y _{k,j} ^{\prime }\colon 2 \le j \le
k)$, the existence
\par\smallskip\noindent
of $\varphi _{k}$ is implied by the embeddability of $K _{k}(Y
_{k,1} ^{\prime })(X _{k,1}, Y _{k,1} ^{\prime }; \varepsilon
_{\mu (k,1)}) _{p ^{\mu (k,k)}}$ as
\par\smallskip\noindent
a $K _{k}$-subalgebra in $K _{k}(X _{k,1}, Y _{k,1}; \varepsilon
_{1+\mu (k,k)}) _{p ^{1+\mu (k,1)}}$. Note also that
\par\smallskip\noindent
$\varphi _{k}(Y _{k,1} ^{\prime }) = \varepsilon _{1+\mu (k,1)}
^{p^{\mu (k,1)}}Y _{k,1} ^{\prime } = \varepsilon _{1}Y _{k,1}
^{\prime }$, and $\varphi _{k} ^{p}(r _{k}) = \xi _{k,1} ^{-1}r
_{k}\xi _{k,1}$, for all $r _{k} \in S _{k,1}$.
\par\smallskip\noindent
Put $X _{\tilde k} ^{\prime } = \sqrt[p]{X _{\tilde k}}$, $Y
_{\tilde k} ^{\prime } = \sqrt[p]{Y _{\tilde k}}$, and consider
the ring $S _{k,1} ^{\prime } := S _{k,1} \langle \theta _{k},
\sigma _{k}, \tau _{k} \rangle $
\par\vskip0.04truecm\noindent
with generators $\theta _{k}, \sigma _{k}$ and $\tau _{k}$ over $S
_{k,1}$, subject to the following relations:
$$\theta _{k} ^{-1}r _{k}\theta _{k} = \varphi _{k}(r _{k}),
\sigma _{k} ^{-1}r _{k}\sigma _{k} = r _{k}, \ {\rm for \ each} \
r _{k} \in S _{k,1}, \theta _{k} ^{p} = \sigma _{k}\xi _{k,1} =
\xi _{k,1}\sigma _{k},$$
$$\theta _{k} ^{p ^{1+\mu (k, 1)}} = X
_{k,1}X _{\tilde k} ^{\prime }, \tau _{k} ^{-1}\theta _{k}\tau
_{k} = \theta _{k}, \tau _{k} ^{-1}\eta _{k,1}\tau _{k} = \eta
_{k,1},$$
$$\tau _{k} ^{-1}\sigma _{k}\tau _{k} = \varepsilon _{\mu (k,
1)}\sigma _{k}, \ {\rm and} \ \tau _{k} ^{p ^{\mu (k, 1)}} = Y
_{\tilde k} ^{\prime }.$$
\par\smallskip\noindent
It is verified by direct calculations that $$\sigma _{k}\eta
_{k,1} = \eta _{k,1}\sigma _{k}, \ \sigma _{k}\theta _{k} = \theta
_{k}\sigma _{k}, \ \sigma _{k} ^{p^{\mu (k,1)}} = X _{\tilde k}
^{\prime }, \ \tau _{k} ^{-1}\xi _{k,1} ^{-1}\tau _{k} =
\varepsilon _{\mu (k,1)}\xi _{k,1} ^{-1},$$
$$\eta _{k,1}\xi _{k,1}\eta _{k,1} ^{-1} = \varepsilon _{\mu (k,
1)}\xi _{k,1}, \eta _{k,1}\xi _{k,1} ^{-1}\eta _{k,1} ^{-1} =
\varepsilon _{\mu (k, 1)} ^{-1}\xi _{k,1} ^{-1},$$
\par\vskip0.09truecm\noindent
and $\eta _{k,1}\tau _{k} ^{-1}.\xi _{k,1} ^{-1} = \xi _{k,1}
^{-1}.\eta _{k,1}\tau _{k} ^{-1}$. These calculations, and the
characterization
\par\vskip0.04truecm\noindent
of an algebra as a tensor product of two subalgebras, given by
Proposition~c
\par\vskip0.04truecm\noindent
of \cite[Sect. 9.2]{P}, indicate that $S _{k,1} ^{\prime } \in d(C
_{\tilde k})$, where $C _{\tilde k}$ denotes the field
\par\smallskip\noindent
$K _{k}(Y _{k,2} ^{\prime }, \dots , Y _{k,k} ^{\prime }, X
_{\tilde k} ^{\prime }, Y _{\tilde k} ^{\prime })$ if $k \ge 2$
(and $C _{\tilde k} = K _{k}(X _{\tilde k} ^{\prime }, Y _{\tilde
k} ^{\prime })$ if $k = 1$). More
\par\vskip0.09truecm\noindent
precisely, one sees that the $C _{\tilde k}$-algebras
$$S _{k,1} ^{\prime } = C _{\tilde k} \langle \xi _{k,1}, \eta
_{k,1}, \theta _{k}, \sigma _{k}, \tau _{k}\rangle = C _{\tilde k}
\langle \eta _{k,1}, \theta _{k}, \sigma _{k}, \tau _{k}\rangle =
C _{\tilde k} \langle \xi _{k,1}, \eta _{k,1}, \theta _{k}, \tau
_{k}\rangle ,$$
$$C _{\tilde k} \langle \theta _{k}, \eta _{k,1}
\rangle \otimes _{C _{\tilde k}} C _{\tilde k} \langle \tau _{k},
\sigma _{k} \rangle \ {\rm and} \ C _{\tilde k} \langle \theta
_{k}, \eta _{k,1}\tau _{k} ^{-1} \rangle \otimes _{C _{\tilde k}}
C _{\tilde k} \langle \xi _{k,1} ^{-1}, \tau _{k} ^{-1} \rangle
$$
\par\smallskip\noindent
are isomorphic. Also, our calculations show that
$$C _{\tilde k} \langle \theta _{k}, \eta _{k,1} \rangle \cong C
_{\tilde k}(X _{k,1}X _{\tilde k} ^{\prime }, Y _{k,1};
\varepsilon _{1+ \mu (k, 1)}) _{p ^{1+\mu (k,1)}},$$
$$C _{\tilde k} \langle \tau _{k}, \sigma _{k} \rangle \cong C
_{\tilde k}(Y _{\tilde k} ^{\prime }, X _{\tilde k} ^{\prime };
\varepsilon _{\mu (k, 1)}) _{p ^{\mu (k, 1)}},$$
$$C _{\tilde k} \langle \theta _{k}, \eta _{k,1}\tau _{k} ^{-1}
\rangle \cong C _{\tilde k}(Y _{k,1}Y _{\tilde k} ^{-1}, (X
_{k,1}X _{\tilde k} ^{\prime }) ^{-1}; \varepsilon _{1+\mu (k,
1)}) _{p ^{1+\mu (k, 1)}}$$ \noindent (observe that
$$C _{\tilde k} \langle
\theta _{k}, \eta _{k,1}\tau _{k} ^{-1} \rangle = C _{\tilde k}
\langle \eta _{k,1}\tau _{k} ^{-1}, \theta _{k} ^{-1}\colon \
\theta _{k} ^{-1}.\eta _{k,1}\tau _{k} ^{-1} = \varepsilon _{1+\mu
(k, 1)}.\eta _{k,1}\tau _{k} ^{-1}.\theta _{k} ^{-1},$$
$$(\eta _{k,1}\tau _{k} ^{-1}) ^{p ^{1+ \mu (k,1)}} = Y _{k,1}Y
_{\tilde k} ^{-1}, \ \theta _{k} ^{-p ^{1+ \mu (k, 1)}} = (X
_{k,1}X _{\tilde k} ^{\prime }) ^{-1} \rangle ),$$
\par\smallskip\noindent
and $C _{\tilde k} \langle \xi _{k,1} ^{-1}, \tau _{k} ^{-1}
\rangle \cong C _{\tilde k}(X _{k,1} ^{-1}, Y _{\tilde k} ^{\prime
-1}; \varepsilon _{\mu (k, 1)}) _{p ^{\mu (k, 1)}}$ as $C _{\tilde
k}$-algebras.
\par
\smallskip
Put $\xi _{\tilde k,k} = \eta _{k,1}\tau _{k} ^{-1}$, $X _{\tilde
k,k} = Y _{k,1}Y _{\tilde k} ^{-1}$, $\eta _{\tilde k,k} = \theta
_{k} ^{-1}$, $Y _{\tilde k,k} = X _{k,1} ^{-p}.X _{\tilde k}
^{-1}$,
\par\smallskip\noindent
$S _{\tilde k,k} = C _{\tilde k} \langle \xi _{\tilde k,k}, \eta
_{\tilde k,k}\rangle $, $\xi _{\tilde k,\tilde k} = \xi _{k,1}
^{-1}$, $X _{\tilde k,\tilde k} = X _{k,1} ^{-1}$, $\eta _{\tilde
k,\tilde k} = \tau _{k} ^{-1}$, $Y _{\tilde k,\tilde k} = Y
_{\tilde k} ^{-1}$,
\par\medskip\noindent
$S _{\tilde k,\tilde k} = C _{\tilde k} \langle \xi _{\tilde
k,\tilde k}, \eta _{\tilde k,\tilde k}\rangle $. It is easily
verified that deg$(S _{\tilde k,k}) = p ^{1+\mu (k, 1)}$ and
\par\smallskip\noindent
deg$(S _{\tilde k,\tilde k}) = p ^{\mu (k, 1)}$; in terms of Lemma
\ref{lemm8.1} (b), one obtains that
\par\medskip\noindent
$\mu (\tilde k, k) = 1 + \mu (k, 1)$, $\mu (\tilde k, \tilde k) =
\mu (k, 1)$. Put $\Sigma _{\tilde k} ^{\prime } = \{X _{\tilde
k,k}, Y _{\tilde k,k}, X _{\tilde k,\tilde k}, Y _{\tilde k,\tilde
k}\}$ if $k = 1$, and in case $k \ge 2$, put $\tilde j = j - 1$,
$\xi _{\tilde k, \tilde j} = \xi _{k,j}$, $\eta _{\tilde k,\tilde
j} = \eta _{k,j}$, $X _{\tilde k,\tilde j} = X _{k,j}$, $Y
_{\tilde k,\tilde j} = Y _{k,j}$, $S _{\tilde k,\tilde j} =
\widetilde S _{k,j} \otimes _{\widetilde C _{k}} C _{\tilde k}$,
for $j = 2, \dots , k$, where
\par\vskip0.09truecm\noindent
$\widetilde C _{k} = K _{0}(X _{k,j}, Y _{k,j} ^{\prime }\colon j
= 2, \dots , k)$, $\Sigma _{\tilde k} ^{\prime } = \{X _{\tilde
k,1}, Y _{\tilde k,1}, \dots , X _{\tilde k,\tilde k}, Y _{\tilde
k,\tilde k}\}$, and
\par\vskip0.09truecm\noindent $\widetilde S _{k,j} =
\widetilde C _{k} \langle \xi _{k,j}, \eta _{k,j}\rangle $, for
each index $j$. It is easy to see that $\Sigma _{\tilde k}
^{\prime }$ is a generating set of the group $\widetilde \Sigma
_{\tilde k}$. Since, by our inductive hypothesis, $K _{k} = K
_{0}(\Sigma _{k} ^{\prime })$ and
\par\vskip0.09truecm\noindent
the set $\Sigma _{k} ^{\prime } = \{X _{k,i}, Y _{k,i}\colon i =
1, \dots , k\}$ is algebraically independent over $K _{0}$,
\par\vskip0.09truecm\noindent
this ensures that $K _{\tilde k} = K _{0}(\Sigma _{\tilde k}
^{\prime })$, $\Sigma _{\tilde k} ^{\prime }$ is algebraically
independent over $K _{0}$ (being a generating set of $K _{\tilde
k}$ over $K _{0}$ of cardinality $2\tilde k = {\rm trd}(K _{\tilde
k}/K _{0})$, see \cite[Ch.~VIII, Theorem~1.1]{L}), and the field
$C _{\tilde k}$ coincides with
\par\vskip0.04truecm\noindent
$K _{0}(X _{\tilde k,i}, Y _{\tilde k,i} ^{\prime }\colon i = 1,
\dots , \tilde k)$, $K _{\tilde k}(Y _{\tilde k,1} ^{\prime },
\dots , Y _{\tilde k, \tilde k} ^{\prime })$, and $Z _{k,1}(X
_{\tilde k} ^{\prime }, Y _{\tilde k} ^{\prime })$, where
\par\vskip0.04truecm\noindent
$Z _{k,1} = K _{k}(Y _{k,j} ^{\prime }\colon j = 2, \dots , k)$ or
$Z _{k,1} = K _{k}$ depending on whether or not $k \ge 2$. Hence,
$K _{\tilde k}/K _{0}$ and $C _{\tilde k}/K _{0}$ are purely
transcendental extensions with trd$(C _{\tilde k}/K _{0}) = {\rm
trd}(K _{\tilde k}/K _{0}) = 2\tilde k$, which allows to deduce
from Kummer
\par\smallskip\noindent
theory that $C _{\tilde k}/K _{\tilde k}$ is a Galois extension
and $\mathcal{G}(C _{\tilde k}/K _{\tilde k}) \cong
(\mathbb{Z}/p\mathbb{Z}) ^{\tilde k}$ (so
\par\smallskip\noindent
$[C _{\tilde k}\colon K _{\tilde k}] = p ^{\tilde k}$). Assuming
now that $R _{\tilde k} = S _{\tilde k,1} \otimes _{C _{\tilde k}}
\otimes _{C _{\tilde k}} \dots \otimes _{C _{\tilde k}} S _{\tilde
k, \tilde k}$, and
\par\smallskip\noindent
putting $Z _{\tilde k} = C _{\tilde k}$, one obtains from
Proposition \ref{prop2.4} that $R _{\tilde k} \in d(Z _{\tilde
k})$ and
\par\vskip0.04truecm\noindent
the $Z _{\tilde k}$-algebra $R _{\tilde k}$ inherits the
properties of the $Z _{k}$-algebra $R _{k}$ described in
\par\vskip0.04truecm\noindent
(8.1). In addition, it follows that $R _{k}$ is isomorphic to a $K
_{k}$-subalgebra of $R _{\tilde k}$, and in the setting of Lemma
\ref{lemm8.1}, $\mu (\tilde k, \tilde j) = \mu (k, j)$, for $j =
2, \dots , k$. Thus it
\par\vskip0.048truecm\noindent
turns out that $\mu (\tilde k, i)\colon i = 1, \dots , \tilde k$,
are positive integers. Moreover, the field
\par\vskip0.048truecm\noindent
$Z _{\tilde k} = Z _{k,1}(X _{\tilde k} ^{\prime }, Y _{\tilde k}
^{\prime })$ is a purely transcendental extension of $Z _{\tilde
k,1}$ with
\par\vskip0.048truecm\noindent
trd$(Z _{\tilde k}/Z _{k,1}) = 2$. This ensures that $Z _{k,1}$
has no algebraic proper extension
\par\vskip0.04truecm\noindent
in $Z _{\tilde k}$, which implies $Z _{\tilde k} \cap Z _{k} = Z
_{k,1}$ (and proves the concluding assertion of Lemma
\ref{lemm8.1} (a), for $n = 2$).
\par\smallskip
Our inductive argument proves the existence of $K _{n}$-algebras
$R _{n}$, $n \in \mathbb{N}$, with the properties required by
Lemma \ref{lemm8.1} (b) and the former part of Lemma \ref{lemm8.1}
(a). It shows that the fields $Z _{n} = Z(R _{n})$, $n \in
\mathbb{N}$, are purely transcendental extensions of $K _{0}$, and
for each $n$, $Z _{n}/K _{n}$ is a Kummer extension admissible by
the former part of Lemma \ref{lemm8.1} (a). Therefore, by
Proposition \ref{prop2.4} and Lemma \ref{lemm8.1} (b), $K _{n}$ is
the only central $K _{n}$-subalgebra of $R _{n}$, for every $n$.
Using repeatedly the same argument, and applying Kummer
\par\vskip0.048truecm\noindent
theory, one obtains step-by-step that $Z _{n+i} \cap K _{n}(Y
_{n,u} ^{\prime }\colon u = 1, \dots , i) = K _{n}$ whenever $n, i
\in \mathbb{N}$ and $i \le n$. This implies $Z _{2n} \cap Z _{n} =
K _{n}\colon n \in \mathbb{N}$, and completes the proof of Lemma
\ref{lemm8.1}.
\end{proof}
\par
\smallskip
We are now in a position to prove Theorem \ref{theo2.1}. Let $K
_{0}$ be an algebraically closed field, $p \in
\mathbb{P}_{K_{0}}$, $K _{\infty }/K _{0}$ a purely transcendental
extension with trd$(K _{\infty }/K _{0}) = \infty $, and $K/K
_{\infty }$ a finite extension. Suppose first that $K = K _{\infty
}$, trd$(K _{\infty }/K _{0})$ is countable, the sets $\Sigma
_{\infty }$ and $\Sigma _{n}\colon n \in \mathbb{N}$, are defined
as in Lemma \ref{lemm8.1}, and $K _{n} = K _{0}(\Sigma _{n})$, for
each $n$. Put $R = \cup _{n=1} ^{\infty } R _{n}$, where $R
_{n}\colon n \in \mathbb{N}$, are taken in agreement with Lemma
\ref {lemm8.1}. It is easily verified that $R$ is a central
division LFD-algebra over $K$ with $[R\colon K]$ countable, and
that finite-dimensional $K$-subalgebras of $R$ are of $p$-power
degrees. We show that $K$ is the unique $K$-subalgebra of $R$
lying in $d(K)$. Assume the opposite, take a $K$-subalgebra $T$ of
$R$ with $T \in d(K)$ and $T \neq K$, and fix a basis $B _{T}$ of
$T$. Clearly, $B _{T} \subset R _{\nu }$, for some $\nu \in
\mathbb{N}$, which can be chosen so that the structural constants
of $T$ with respect to $B _{T}$ lie in $K _{\nu }$. Then $R _{\nu
}$ has a $K _{\nu }$-subalgebra $T _{\nu }$, such that $T _{\nu }
\otimes _{K _{\nu }} K \cong T$ as $K$-subalgebras. As $T \neq K$,
$T \in d(K)$ and $[T _{\nu }\colon K _{\nu }] = [T\colon K]$, this
requires $T _{\nu } \in d(K _{\nu })$ and $T _{\nu } \neq K _{\nu
}$, which contradicts Proposition \ref{prop2.4} and thereby proves
Theorem \ref{theo2.1} in the case of $K = K _{\infty }$ and trd$(K
_{\infty }/K _{0})$ countable.
\par\smallskip
Assume now that $K \neq K _{\infty }$ or trd$(K _{\infty }/K
_{0})$ is uncountable and fix an algebraic closure $\overline K$
of $K$. Then $K _{0}$ has purely transcendental extensions $\Psi $
and $\Phi $ in $K _{\infty }$, such that trd$(\Psi /K _{0})$ is
countable, $\Psi \cap \Phi = K _{0}$,
\par\vskip0.04truecm\noindent
$\Psi \Phi = K _{\infty } \cong \Psi \otimes _{K _{0}} \Phi $, and
there exists $\Phi ^{\prime } \in I(K/\Phi )$ satisfying
\par\vskip0.04truecm\noindent
$[\Phi ^{\prime }\colon \Phi ] = [K\colon K _{\infty }]$ and $\Phi
^{\prime }\Psi = K$. This ensures that $\Phi ^{\prime } \in
I(\overline \Phi /\Phi )$,
\par\vskip0.04truecm\noindent
$\overline \Phi \Psi \cong \overline \Phi \otimes _{K _{0}} \Psi $
and $\overline \Phi \Psi /\overline \Phi $ is a purely
transcendental extension with
\par\vskip0.04truecm\noindent
trd$(\overline \Phi \Psi /\overline \Phi ) = {\rm trd}(\Psi /K
_{0})$, where $\overline \Phi $ is the algebraic closure of $\Phi
$ in $\overline K$. Moreover, every system of elements of $\Psi $
algebraically independent over $K _{0}$ remains algebraically
independent over $\overline \Phi $. Therefore, applying Lemma
\ref{lemm8.1} and arguing as in the preceding part of our proof,
one proves the existence of a central division LFD-algebra $R
_{\Psi }$ over $\Psi $, such that $[R _{\Psi }\colon \Psi ]$ is
countable, $R _{\Psi }$ does not possess noncommutative $\Psi
$-subalgebras lying in $d(\Psi )$, and $R _{\Psi } \otimes _{\Psi
} \overline \Phi \Psi $ preserves the same properties as a
$\overline \Phi \Psi $-algebra. Let now $\Phi _{1}$ be an
arbitrary extension of $\Phi $ in $\overline \Phi $. Then
\par\vskip0.04truecm\noindent
$R _{\Psi } \otimes _{\Psi } \overline \Phi \Psi \cong (R _{\Psi }
\otimes _{\Psi } \Phi _{1}\Psi ) \otimes _{\Phi _{1}\Psi }
\overline \Phi \Psi $ as $\overline \Phi \Psi $-algebras (see
\cite[Sect.~9.4, Corollary~a]{P}), so it can be easily obtained
from the noted properties of
\par\vskip0.05truecm\noindent
$R _{\Psi } \otimes _{\Psi } \overline \Phi \Psi $ that $R _{\Psi
} \otimes _{\Psi } \Phi _{1}\Psi $ is a central division
LFD-algebra over $\Phi _{1}\Psi $,
\par\vskip0.048truecm\noindent
$[(R _{\Psi } \otimes _{\Psi } \Phi _{1}\Psi )\colon \Phi _{1}\Psi
]$ is countable, and the class $d(\Phi _{1}\Psi )$ does not
contain
\par\vskip0.048truecm\noindent
noncommutative $\Phi _{1}\Psi $-subalgebras of $R _{\Psi } \otimes
_{\Psi } \Phi _{1}\Psi $. This, applied to  the case
\par\vskip0.04truecm\noindent
of $\Phi _{1} = \Phi ^{\prime }$ (and $\Phi ^{\prime }\Psi = K$),
completes the proof of Theorem \ref{theo2.1}.
\par
\medskip
We conclude with an analog to Proposition \ref{prop3.5} for tensor
products of central division LFD-algebras admissible by Theorem
\ref{theo2.1} (a) and (b).
\par
\medskip
\begin{prop}
\label{prop8.2} With assumptions being as in Proposition
\ref{prop3.5}, suppose that $d(K)$ does not contain a
noncommutative $K$-subalgebra of $R\{p\}$, for
\par\noindent
any $p \in \mathbb{P}$. Then $R = \otimes _{p \in \mathbb{P}}
R\{p\}$ is a central division {\rm LFD}-algebra over $K$ whose
noncommutative $K$-subalgebras do not belong to $d(K)$.
\end{prop}
\par
\smallskip
\begin{proof}
We prove that if $d(K)$ does not contain any noncommutative
$K$-subalgebra of $R\{p\}$, for a given $p \in \mathbb{P}$, then
$p \nmid {\rm deg}(\mathcal{R}_{p})$, for any $K$-subalgebra
$\mathcal{R}_{p}$ of $R$ with $\mathcal{R}_{p} \in d(K)$
(Proposition \ref{prop8.2} follows at once from this statement).
Using Brauer's theorem, one observes that the statement will be
proved if we show that a $K$-subalgebra $\Sigma _{p}$ of $R$
embeds in the $K$-algebra $R\{p\}$, provided that $\Sigma _{p} \in
d(K)$ and deg$(\Sigma _{p}) = p ^{\sigma }$, for some $\sigma \in
\mathbb{N}$.
\par\vskip0.04truecm\noindent
Let $R ^{\prime }$ and $R\{p\} ^{\prime }$ be the underlying
division $K$-algebras of $R \otimes _{K} \Sigma _{p} ^{\rm op}$
and
\par\vskip0.04truecm\noindent
$R\{p\} \otimes _{K} \Sigma _{p} ^{\rm op}$, respectively
(determined by the Wedderburn-Artin theorem). Then it is not
difficult to see that $R ^{\prime }$ and $R\{p\} ^{\prime }$ are
central LFD-algebras over $K$. In addition, Lemma~3.5 of
\cite{Ch2} and the embeddability of $\Sigma _{p}$ in $R$ imply
that $R \otimes _{K} \Sigma _{p} ^{\rm op}$ and the matrix ring $M
_{p^{2\sigma }}(R ^{\prime })$ are isomorphic as $K$-algebras;
also, by the same lemma, there is a $K$-isomorphism
\par\vskip0.04truecm\noindent
$R\{p\} \otimes _{K} \Sigma _{p} ^{\rm op} \cong M_{p^{y}}(R\{p\}
^{\prime })$, for some integer $y \ge 0$. We show that every
finite-dimensional $K$-subalgebra $\Lambda _{p} ^{\prime }$ of
$R\{p\} ^{\prime }$ has $p$-power dimension. Relying upon the
basic theory of finite-dimensional central division algebras (cf.
\cite[Sects.~12.1 and 13.1]{P}), one may consider only the special
case where $\Lambda _{p} ^{\prime }$ is a field. Note further that
one may assume, for the proof, that $\Lambda _{p} ^{\prime }/K$ is
a separable extension. This is evident if char$(K) = 0$. When $p =
{\rm char}(K)$, the reduction is a result of the fact that $p \mid
[K _{1}\colon K]$, for any inseparable finite extension $K
_{1}/K$, and in the case where $p \neq {\rm char}(K) \neq 0$, it
is obtained by the method of proving the Noether-Jacobson theorem
(e.g., in \cite[Ch.~3]{He}). Let now $L$ be any finite extension
of $K$ in $K _{\rm sep}$, $M$ the Galois closure of $L$ in $K
_{\rm sep}$ over $K$, $P ^{\prime }$ a Sylow $p$-subgroup of
$\mathcal{G}(M/L)$, $P$ a Sylow $p$-subgroup of $\mathcal{G}(M/K)$
including $P ^{\prime }$, $F ^{\prime }$ and $F$ the fixed fields
of $P ^{\prime }$ and $P$, respectively. It is clear from Galois
theory and Sylow's theorems that $p \nmid [F\colon K]$. This,
combined with Lemma \ref{lemm3.6}, proves that $\Lambda _{p}
\otimes _{K} F$ is a division $F$-algebra, for every
finite-dimensional $K$-subalgebra $\Lambda _{p}$ of $R\{p\}$; thus
$R\{p\} \otimes _{K} F$ turns out to be a central division
$F$-algebra. Similarly, one concludes that $\Sigma _{p} ^{\rm op}
\otimes _{K} F \in d(F)$, and then deduces from
\cite[Lemma~3.5]{Ch2} that there exists an $F$-isomorphism
$(R\{p\} \otimes _{K} F) \otimes _{F} (\Sigma _{p} ^{\rm op}
\otimes _{K} F) \cong M_{p^{\tau }}(\widetilde R\{p\})$, for some
integer $\tau \ge 0$ and a central division $F$-algebra
$\widetilde R\{p\}$. At the same time, by \cite[Lemma~3.5]{Ch2},
$R\{p\} ^{\prime } \otimes _{K} F$ is isomorphic as an $F$-algebra
to $M _{\gamma }(\widetilde R\{p\} ^{\prime })$, for some $\gamma
\in \mathbb{N}$ dividing $[F\colon K]$, and some central division
$F$-algebra $\widetilde R\{p\} ^{\prime }$. Therefore, it is not
difficult to show that there are $F$-isomorphisms
$$(R\{p\} \otimes _{K} \Sigma _{p}
^{\rm op}) \otimes _{K} F \cong (R\{p\} \otimes _{K} F) \otimes
_{F} (\Sigma _{p} ^{\rm op} \otimes _{K} F) \cong M _{p^{\tau
}}(\widetilde R\{p\})$$ $$\cong M_{p^{y}}(R\{p\} ^{\prime })
\otimes _{K} F \cong M _{p^{y}}(F) \otimes _{F} M_{\gamma
}(\widetilde R\{p\} ^{\prime }) \cong M_{p^{y}.\gamma }(\widetilde
R\{p\} ^{\prime }).$$ In view of the Wedderburn-Artin theorem,
this requires that $y = \tau $, $\gamma = 1$
\par\smallskip\noindent
and $\widetilde R\{p\} \cong \widetilde R\{p\} ^{\prime }$ as
$F$-algebras, proving that $R\{p\} ^{\prime } \otimes _{K} F \cong
\widetilde R\{p\}$. Since
\par\smallskip\noindent
$L \in I(F ^{\prime }/K)$ and $[F ^{\prime }\colon F] = p ^{z}$,
for some integer $z \ge 0$, and there exist
\par\medskip\noindent
$F ^{\prime }$-isomorphisms $R\{p\} ^{\prime } \otimes _{K} F
^{\prime } \cong (R\{p\} ^{\prime } \otimes _{K} F) \otimes _{F} F
^{\prime } \cong (R\{p\} ^{\prime } \otimes _{K} L) \otimes _{L} F
^{\prime }$
\par\medskip\noindent
(cf. \cite[Sect.~9.4, Corollary~a]{P}), one also sees that $R\{p\}
^{\prime } \otimes _{K} L \cong M_{p^{t}}(\mathcal{R}\{p\})$
\par\medskip\noindent
and $R\{p\} ^{\prime } \otimes _{K} F ^{\prime } \cong
M_{p^{t'}}(\mathcal{R}\{p\} ^{\prime })$, for some integers $t$
and $t'$ with
\par\smallskip\noindent
$0 \le t \le t' \le z$. Here $\mathcal{R}\{p\}$ is the underlying
central division $L$-algebra
\par\smallskip\noindent
of $R\{p\} ^{\prime } \otimes _{K} L$, and $\mathcal{R}\{p\}
^{\prime }$ is the underlying central division $F ^{\prime
}$-algebra of $\mathcal{R}\{p\} \otimes _{L} F ^{\prime }$,
determined by the Wedderburn-Artin theorem. When $L$ is a
\par\vskip0.04truecm\noindent
$K$-isomorphic copy of $\Lambda _{p} ^{\prime }$, these
calculations and \cite[Lemma~3.5]{Ch2} prove
\par\smallskip\noindent
that $[\Lambda _{p} ^{\prime }\colon K]$ is a $p$-power, as
claimed. As $Z(R\{p\} ^{\prime }) = Z(R\{p\}) = K$ and $R\{p\}
^{\prime }$ is LFD over $K$, Proposition \ref{prop3.5} and this
fact indicate that
\par\smallskip\noindent
$R\{p\} ^{\prime } \otimes _{K} (\otimes _{p' \in \mathbb{P} _{p}}
R\{p ^{\prime }\})$ is a central division LFD-algebra over $K$,
where $\mathbb{P} _{p} = \mathbb{P} \setminus \{p\}$. The obtained
result and the $K$-isomorphism
\par\smallskip\noindent
$R \otimes _{K} \Sigma _{p} ^{\rm op} \cong M _{p^{2\sigma }}(R
^{\prime })$ enable one to deduce from the Wedderburn-Artin
\par\smallskip\noindent
theorem that $R\{p\} ^{\prime } \otimes _{K} (\otimes _{p' \in
\mathbb{P} _{p}} R\{p ^{\prime }\}) \cong R ^{\prime }$ and
$R\{p\} \otimes _{K} \Sigma _{p} ^{\rm op} \cong M _{p^{2\sigma
}}(R\{p\} ^{\prime })$
\par\smallskip\noindent
as $K$-algebras. It is now easy to conclude that $y = 2\sigma $,
and by \cite[Lemma~3.5]{Ch2}, $\Sigma _{p}$ is embeddable in
$R\{p\}$ as a $K$-subalgebra, which proves Proposition
\ref{prop8.2}.
\end{proof}
\par
\medskip
It is unknown whether every central division LFD-algebra $R$ over
a field $K \notin \Phi _{\rm Br}$ with $[R\colon K]$ countable
admits a primary tensor product decomposition. Since, by
K\"{o}the's theorem (generalized by \cite[Theorem~13.22]{BGMV}),
the answer is affirmative if $R$ is NLF over $K$, a negative
answer in general would lead to new solutions to the normality
problem with different structure from those in the present paper.
When $[R\colon K]$ is uncountable (which may occur if Brd$_{p}(K)
\neq 0$, for infinitely many $p$), a negative answer is also
possible if $R$ is NLF, see \cite[(5)]{Ba}. The same holds for any
$K$ and other classes of central simple LFD-algebras, such as
locally matrix $K$-algebras (see \cite[Theorem~2]{Ku} and
\cite[Theorem~2]{BO}).
\par\vskip0.14truecm
\emph{Acknowledgement.} This research has partially been supported
by the Bulgarian National Science Fund, Grant KP-06 N 32/1 of
07.12.2019, and by the Science Foundation of Sofia University
under contract 80-10-68/9.4.2024.
\par\vskip0.14truecm

\par

\end{document}